\documentclass[english,11pt,leqno,twoside]{amsart}
\usepackage[T1]{fontenc}
\usepackage[utf8]{inputenc}
\usepackage{amsmath}
\usepackage{amsthm}
\usepackage{amssymb}

\usepackage{soul,color}
\usepackage{t1enc}
\usepackage{a4,indentfirst,latexsym,graphics}
\usepackage{mathrsfs}
\usepackage[noadjust]{cite}
\usepackage{enumitem,graphicx}
\usepackage[colorlinks=true,urlcolor=blue,
citecolor=red,linkcolor=blue,linktocpage,pdfpagelabels,
bookmarksnumbered,bookmarksopen]{hyperref}
\usepackage[english]{babel}
\usepackage[left=2.1cm,right=2.1cm,top=2.72cm,bottom=2.72cm]{geometry}
\usepackage[metapost]{mfpic}
\usepackage[hyperpageref]{backref}
\usepackage[colorinlistoftodos]{todonotes}
\usepackage[normalem]{ulem}
\usepackage{braket}
\usepackage{fancyhdr}
\makeatletter
\providecommand\@dotsep{5}
\def\listtodoname{List of Todos}
\def\listoftodos{\@starttoc{tdo}\listtodoname}
\makeatother

\newtheorem{Th}{Theorem}[section]
\newtheorem{Prop}[Th]{Proposition}
\newtheorem{Lem}[Th]{Lemma}
\newtheorem{Cor}[Th]{Corollary}
\theoremstyle{definition}
\newtheorem{Def}[Th]{Definition}

\theoremstyle{remark}
\newtheorem{Rem}[Th]{Remark}

\newcommand{\R}{\mathbb{R}}
\newcommand{\rn}{\R^N}
\newcommand{\hrn}{H^1(\rn)}

\newcommand{\cS}{\mathcal{S}}
\newcommand{\cM}{\mathcal{M}}
\newcommand{\cD}{\mathcal{D}}
\newcommand{\cC}{\mathcal{C}}

\newcommand{\cG}{\mathcal{G}}
\newcommand{\cU}{\mathcal{U}}
\newcommand{\dx}{\mathrm{d}x}
\newcommand{\dt}{\mathrm{d}t}
\newcommand{\ds}{\mathrm{d}s}

\newcommand{\tu}{\widetilde{u}}
\newcommand{\hrnr}{H^1_{\mathrm{rad}}(\rn)}

\DeclareMathOperator*{\esssup}{ess\,sup}

\definecolor{yellow-green}{rgb}{0.6, 0.8, 0.2}

\numberwithin{equation}{section}

\title[Normalized solutions to Schr\"odinger equations with double-behaviour nonlinearities]{Existence and dynamics of normalized solutions to Schr\"odinger equations with generic double-behaviour nonlinearities}

\author[Bartosz Bieganowski]{Bartosz Bieganowski}
\author[Pietro d'Avenia]{Pietro d'Avenia}
\author[Jacopo Schino]{Jacopo Schino}

\address[B. Bieganowski]{\newline\indent
	Faculty of Mathematics, Informatics and Mechanics, \newline\indent
	University of Warsaw, \newline\indent
	ul. Banacha 2, 02-097 Warsaw, Poland}
\email{\href{mailto:bartoszb@mimuw.edu.pl}{bartoszb@mimuw.edu.pl}}

\address[P. d'Avenia]{\newline\indent
	Dipartimento di Meccanica, Matematica e Management, \newline\indent
	Politecnico di Bari, \newline\indent
	Via E. Orabona 4, 70125 Bari, Italy}
\email{\href{mailto:pietro.davenia@poliba.it}{pietro.davenia@poliba.it}}

\address[J. Schino]{\newline\indent
	Faculty of Mathematics, Informatics and Mechanics, \newline\indent
	University of Warsaw, \newline\indent
	ul. Banacha 2, 02-097 Warsaw, Poland}
\email{\href{mailto:j.schino2@uw.edu.pl}{j.schino2@uw.edu.pl}}

\subjclass[2020]{35Q55, 35Q40, 35J20}

\keywords{
	Nonlinear Schr\"odinger equations,
	normalized solutions,
	least-energy solutions,
	multiple solutions, 
	mixed nonlinearities,
	variational methods,
	orbital stability,
	strong instability}

\begin{document}
	
	\begin{abstract}
		We study the existence of solutions $(\underline u,\lambda_{\underline u})\in H^1(\mathbb{R}^N; \mathbb{R}) \times \mathbb{R}$ to 
		\[
		-\Delta u + \lambda u = f(u) \quad \text{in } \mathbb{R}^N
		\]
		with $N \ge 3$ and prescribed $L^2$ norm, and the dynamics of the solutions to
		\[
		\begin{cases}
			\mathrm{i} \partial_t \Psi + \Delta \Psi = f(\Psi)\\
			\Psi(\cdot,0) = \psi_0 \in H^1(\mathbb{R}^N; \mathbb{C})
		\end{cases}
		\]
		with $\psi_0$ close to $\underline u$.
		Here, the nonlinear term $f$ has mass-subcritical growth at the origin, mass-supercritical growth at infinity, and is more general than the sum of two powers. Under different assumptions, we prove the existence of a locally least-energy solution, the orbital stability of all such solutions, the existence of a second solution with higher energy, and the strong instability of such a solution.
	\end{abstract}
	
	\maketitle

	\begin{center}
		\begin{minipage}{11cm}
			\tableofcontents
		\end{minipage}
	\end{center}

	\section{Introduction}

	The nonlinear Schr\"odinger equation with a nonlinearity consisting of combined powers, i.e.,
	\[
	\mathrm{i} \partial_t \Psi + \Delta \Psi = \mu_1 |\Psi|^{p_1-2} \Psi + \mu_2 |\Psi|^{p_2-2} \Psi, \quad \Psi \colon \R \times \R^N \to \mathbb{C},
	\]
	with $\mu_1,\mu_2 \in \R \setminus \{0\}$ and $2 < p_1 < p_2 \le \frac{2N}{N-2}$, attracted much attention since when the pioneering paper by T.C. Tao, M. Vi\c{s}an, and X. Zhang \cite{TVZ} appeared. After that, several results have been obtained in the context of solutions with prescribed $L^2$-norm. 
	
	The seminal work was done by Soave, who studied the existence and properties of solutions in the Sobolev-subcritical \cite{Soave} and Sobolev-critical \cite{SoaveC} cases, i.e.,
	\begin{equation*}
		\begin{cases}
			-\Delta u + \lambda u = \mu |u|^{q-2} u + |u|^{p-2} u,\\
			\displaystyle \int_{\R^N} |u|^2 \, \dx = \rho^2,\\
			(u,\lambda) \in H^1(\R^N;\mathbb{R}) \times \R
		\end{cases}
	\end{equation*}
	with $\rho>0$, $\mu \in \R \setminus \{0\}$, $2 < q \le 2_\# \le p$, $q < p$, and either $p < 2^*$ or, when $N \ge 3$, $p = 2^*$, where $2^*:=2N/(N-2)$ is the Sobolev-critical exponent
	and
	\[
	2_\# := 2 + \frac4N
	\]
	is the so-called {\em mass-} or {\em $L^2$-critical exponent}. We point out that $2_\#$ often plays a crucial role in the dynamics of the nonlinear Schr\"odinger equation and in the context of normalized solutions.
	
	Some questions in \cite{SoaveC} concerning the existence and asymptotic behaviour of mountain-pass solutions remained open and were answered later on in \cite{JJLV,Jeanjean_Le,Qi_Zou,Wei_Wu}.

	A common feature of \cite{Jeanjean_Le,Qi_Zou,Soave,SoaveC,Wei_Wu} is that the elliptic problem is proved to have a least-energy solution and, in some cases, a second solution of mountain-pass-type. In the recent paper \cite{AJM}, however, the authors obtain an arbitrarily large number of solutions with prescribed $L^2$-norm in the presence of Sobolev-critical ($N \ge 3$) or exponential critical ($N=2$) nonlinear terms.
	
	Similar problems (i.e., where the nonlinear term behaves differently at the origin and at infinity with respect to the mass-critical exponent) were studied as well:
	we mention, e.g., systems of equations \cite{Bartsch_Jeanjean,BLZ,Gou_Jeanjean,hLi_Zou_JFPTA,LWYZ}, non-autonomous equations \cite{BQZ,Gou,hLi_Zou_MN,LXZ}, or equations on domains \cite{BQZ,Qi_Zou_bbd}. We would also like to bring \cite{Jeanjean_Lu} to attention, which considers nonlinearities $f(t)$ with different behaviours at the origin and at infinity, but in a situation opposite to \cite{Soave,SoaveC}: $\limsup_{|t| \to +\infty} f(t)t / |t|^{2_\#} \le 0$ and $\limsup_{t \to 0} f(t)t / |t|^{2_\#} \le 0$.
	
	Concerning the non-mixed case, the interested reader will find plenty of results in the literature. Here, we simply mention some pioneering work: \cite{Lions,Stuart82} in the case where the energy -- see \eqref{eq:energy} -- is bounded below, and \cite{Jeanjean} in the case where it is not.
	
	Finally, we report \cite{JZZ}, which deals with all the possible regimes (for an autonomous equation) via a bifurcation method.
	
	Inspired by this work, we are interested in the following normalized problem
	\begin{equation}\label{eq:main}
		\begin{cases}
			-\Delta u + \lambda u = f(u),\\
			\displaystyle \int_{\R^N} |u|^2 \,\dx = \rho^2,\\
			(u,\lambda) \in H^1(\R^N;\R) \times \R,
		\end{cases}
	\end{equation}
	where $\rho>0$ is prescribed and $f \colon \R \rightarrow \R$ is a general nonlinear function that behaves in a way similar to the sum of two powers. In particular, $f$ is mass-subcritical at $0$ and (possibly) mass-supercritical at infinity, cf. assumptions (\ref{F0})--(\ref{F4}), (\ref{H0n}), and (\ref{H2n}).

	Since we do not (always) assume that $f$ is odd, we work in $H^1(\rn) = H^1(\rn;\R)$. However, if we assume that $f$ is continuous and odd, then, as usual (cf. \cite[Example 3.2.4]{Cazenave}), we can easily extend it to the complex plane setting 
	\begin{equation*}
		f(z)
		=
		\begin{cases}
			f(|z|)z/|z| & \text{for } z \in \mathbb{C} \setminus \{0\}\\
			0 & \text{for } z=0.
		\end{cases}
	\end{equation*}
	In this case, we can work in $H^1(\R^N;\mathbb{C})$. Conversely, whenever we mention functions in $H^1(\R^N;\mathbb{C})$, we implicitly assume that the real-valued $f$ is odd and so extended to the complex plane as above. With this framework, \eqref{eq:main} is closely related to the Cauchy problem
	\begin{equation}\label{eq:time}
		\begin{cases}
			\mathrm{i} \partial_t \Psi + \Delta \Psi = f(\Psi)\\
			\Psi(\cdot,0) = \psi_0 \in H^1(\rn; \mathbb{C}).
		\end{cases}
	\end{equation}
	Indeed, a standing wave $\Psi(t,x) = u(x) e^{\mathrm{i}\lambda t}$ satisfies \eqref{eq:time} if and only if $(u,\lambda)$ is a solution to \eqref{eq:main} with $\rho = |\psi_0|_2 = |u|_2$, where $|\cdot|_k$ denotes the usual $L^k$-norm.
	
	Problems like this have applications in mathematical physics. They are used to describe various physical phenomena, including quantum mechanics and fluid dynamics. In quantum mechanics, nonlinear equations similar to \eqref{eq:time} are used to understand the behaviour of subatomic particles and design quantum technologies \cite{ReedSimon}. Other applications include nonlinear optics \cite{optics} and Bose--Einstein condensates \cite{BEcond}.

	We introduce the following assumptions about $f$, where, here and in what follows, $\lesssim$ denotes an inequality up to a positive multiplicative constant.

	\begin{enumerate}[label=(F\arabic{*}),ref=F\arabic{*}]\setcounter{enumi}{-1}
		\item \label{F0} $f$ is continuous and $|f(t)| \lesssim |t| + |t|^{2^*-1}$.
		\item \label{F1} $\displaystyle \lim_{t \to 0} \frac{F(t)}{t^2} = 0$.
		\item \label{F2} $\displaystyle \lim_{t \to 0} \frac{F(t)}{|t|^{2_\#}} = +\infty$.
		\item \label{F3} $\displaystyle \lim_{|t| \to +\infty} \frac{F(t)}{|t|^{2^*}} = 0$.
		\item \label{F4} $f(t)t \leq 2^* F(t)$ for all $t \in \mathbb{R}$.
	\end{enumerate}
	
	\begin{Rem}
		From (\ref{F0}), if $F(\zeta) > 0$ for some $\zeta \neq 0$ (which occurs if (\ref{F2}) holds), the number
		\begin{equation}\label{eq:C0}
			C_0 := \sup_{0 \ne t \in \R} \frac{F(t)}{t^2 + |t|^{2^*}} > 0
		\end{equation}
		is well-defined. We point out that, hereafter, whenever we mention $C_0$, we implicitly assume that $F$ is positive somewhere and $C_0 <+\infty$. 
	\end{Rem}
	
	The \textit{energy functional} $J\colon \hrn \to \R$ associated with \eqref{eq:main} is given by
	\begin{equation}\label{eq:energy}
		J(u) = \frac12 \int_{\rn} |\nabla u|^2 \, \dx - \int_{\rn} F(u) \, \dx
	\end{equation}
	with $F(t) = \displaystyle \int_0^t f(s) \, \mathrm{d}s$, where in $H^1(\R^N)$ we consider the usual norm $\|u\| := ( |\nabla u|_2^2 + |u|_2^2 )^{1/2}$.\\
	If $f$ is odd, then $F$ is even, so $J$ naturally extends to the energy $J \colon H^1 (\R^N; \mathbb{C}) \rightarrow \R$ associated with \eqref{eq:time} by setting $F(z) := F(|z|)$, $z \in \mathbb{C}$. Note that, under mild assumptions -- cf. (\ref{F0}), the mass and the energy of solutions to \eqref{eq:time} are conserved in time, i.e.,
	$$
	|\Psi(\cdot, t)|_2 = |\psi_0|_2, \quad J(\Psi(\cdot, t)) = J(\psi_0).
	$$
	Moreover, by (\ref{F4}) we have $0 \leq f(z) \overline{z} \leq 2^* F(z)$ for all $z \in \mathbb{C}$.

	For $\alpha > 0$, let
	\begin{equation*}
		\cD_\alpha := \Set{u \in \hrn | \left|u \right|_2 \le \alpha} \quad \text{and} \quad
		\cS_\alpha := \Set{u \in \hrn | \left|u\right|_2 = \alpha}.
	\end{equation*}
	Additionally, for $\alpha,R > 0$, let us introduce
	\begin{equation}\label{URamRa}
		\cU_R(\alpha) := \Set{u \in \cD_\alpha | \left|\nabla u\right|_2 < R} \quad \text{and} \quad m_R(\alpha) := \inf_{\cU_R(\alpha)} J.
	\end{equation}
	
	The idea of working with $\cD_\alpha$ instead of $\cS_\alpha$ was introduced in \cite{BiegMed} in the context of nonlinearities with mass-critical or -supercritical growth at the origin and mass-supercritical and Sobolev-subcritical at infinity. Later on, it was exploited in different scenarios, see \cite{BMS,CLY,Liu_Zhao_CVPD,MedSch,MSlog,Schino,LRZ}. The main advantage is that the weak limit of a sequence in $\cD_\alpha$ still belongs to $\cD_\alpha$, while this is not the case with $\cS_\alpha$ because the embedding $\hrn \hookrightarrow L^2(\R^N)$ is not compact, not even considering radially symmetric functions. This makes it easier to obtain a minimizer of $J$ over suitable subsets, which is an important step to obtain a solution to \eqref{eq:main}.
	
	To state our results, we recall the Gagliardo--Nirenberg--Sobolev inequality
	\begin{equation}\label{eq:gns}
		S|u|_{2^*}^{2} \leq |\nabla u|_2^{2}, \quad u \in D^{1,2}(\R^N),
	\end{equation}
	where $S := \inf\set{\left|\nabla u\right|_2^{2}/\left|u\right|_{2^*}^{2} | u \in D^{1,2}(\R^N) \setminus \{0\}} > 0$ is the best constant and $D^{1,2}(\rn)$ denotes the completion of $\cC_0^\infty(\rn)$ with respect to the norm $|\nabla \cdot|_2$.
	
	Our first result is the following.
	\begin{Th}\label{th:locmin}
		Assume (\ref{F0})--(\ref{F2}) and
		\begin{equation}\label{eq:rho}
			\rho^2 < \frac{2}{N-2} \left(\frac{S}{2^* C_0}\right)^{N/2}.
		\end{equation}
		If $\rho_n \to \rho$, $R_n \to R_0$ (where $R_0 > 0$ is defined in (\ref{eq:g_prop})), and $u_n \in \cU_{R_n}(\rho_n)$ are such that $J(u_n) \to m_{R_0}(\rho)$, then there exist $\overline u \in \cS_\rho\cap\cU_{R_0}(\rho)$ and $\lambda_{\overline u} > 0$ such that $u_n \to \overline u$ up to translations and up to subsequences, $J(\overline u) = m_{R_0}(\rho) < 0$, $\overline{u}$ has constant sign, is nowhere zero, and $(\overline u,\lambda_{\overline u})$ is a solution to \eqref{eq:main}.\\
		Actually, if $u \in \cU_{R_0}(\rho)$ satisfies $J(u) = m_{R_0}(\rho)$, then $u$ has constant sign, is nowhere zero, and $(u,\lambda_u)$ is a solution to \eqref{eq:main} for some $\lambda_u > 0$.
	\end{Th}

	\begin{Rem}\label{rem:12}
		With the assumptions and notations of Theorem \ref{th:locmin}, the following facts hold.
		\begin{enumerate}[label=(\roman{*}),ref=\roman{*}]
			\item \label{Linfty} From a standard elliptic argument, $u \in L^\infty(\R^N) \cap \cC^{1,\alpha}(\R^N)$ for every $\alpha \in (0,1)$ and $u(x) \to 0$ as $|x| \to +\infty$. We will need this in some of the proofs below.
			\item One can always assume that $u$ is radial and radially monotonic (in particular, non-increasing if non-negative, non-decreasing if non-positive) because, since it has constant sign, one can replace it with its Schwarz rearrangement \cite[Chapter 3]{LiebLoss} (if $u \ge 0$) or the negative of its Schwarz rearrangement (if $u \le 0$).
			\item If $u$ is radial and radially monotonic and there exist $\tau_+,\tau_- > 0$ such that $f(t) \le \lambda_u t$ for all $t \in [0,\tau_+]$, $f(t) > \lambda_u t$ for all $t > \tau_+$, $f(t) \ge \lambda_u t$ for all $t \in [-\tau_-,0]$, and $f(t) < \lambda_u t$ for all $t < -\tau_-$, then $u$ is radially strictly monotonic; cf. the proof of Theorem \ref{th:2sol} below.
		\end{enumerate}
	\end{Rem}

	Somewhat surprisingly, the only assumption we need for Theorem \ref{th:locmin} (other than those creating the geometry of $J$, which defines the regime we study here) is that $\rho$ is sufficiently small. In particular, we do not distinguish between nonlinear terms that have Sobolev-critical or -subcritical growth at infinity. In fact, one of the purposes of this paper is to understand what reasonably minimal hypotheses we need both for the various steps (this will be evident especially in Section \ref{sec:pes}) and for the main results.
	
	We introduce the following abstract assumptions, where
	\begin{equation}\label{star}
		s \star u := s^{N/2} u(s\cdot)
	\end{equation}
	for $u \in \hrn \setminus \{0\}$ and $s>0$.
	\begin{enumerate}[label=(J\arabic{*}),ref=J\arabic{*}]
		\item \label{J1} For every $u \in \cD_\rho \setminus \{0\}$ the function $(0,+\infty) \ni s \mapsto J(s\star u) \in \R$ has a unique local maximum point $t_u$.
		\item \label{J2} Assumption (\ref{J1}) holds and for every $u \in \cD_\rho \setminus \{0\}$ the function $(t_u,+\infty) \ni s \mapsto J(s\star u) \in \R$ is concave.
	\end{enumerate}
	An example that satisfies (\ref{J2}) is $f(t) = \nu |t|^{p-2}t + \mu |t|^{q-2}t$ with $\mu,\nu > 0$ and $2 < q < 2_\# < p \le 2^*$, cf. \cite{Soave,SoaveC}. Observe that (\ref{J2}) also holds in a more general mass-supercritical regime under some Ambrosetti--Rabinowitz-type assumptions, see \cite{BiegMed}. Other examples satisfying (\ref{J1}) or (\ref{J2}) are shown in Appendix \ref{Example}; in particular, the value $2 + 2/N$ plays a critical role and appears to be a new, important threshold.
	
	It makes sense to ask whether the local minimizer found in Theorem \ref{th:locmin} is a least-energy solution. This is indeed the case under the additional assumption (\ref{J2}).
	
	\begin{Prop}\label{pr:gs}
		If (\ref{F0})--(\ref{F2}), (\ref{J2}), and \eqref{eq:rho} hold, then
		\begin{equation*}
			m_{R_0}(\rho) = \min \Set{J(u) | u \in \cD_\rho \text{ and } J|_{\cD_\rho}'(u) = 0}.
		\end{equation*}
	\end{Prop}
	
	\begin{Rem}
		Note that if $u \in \cS_\rho$, then $J|_{\cS_\rho}'(u) = J'(u) + \lambda_u u$ for some $\lambda_u \in \R$. Hence, above and in what follows, we understand $J|_{\cD_\rho}'(u)$ as
		\begin{align*}
			J|_{\cD_\rho}'(u) := \left\{ \begin{array}{ll}
				J|_{\cS_\rho}'(u) & \quad \mbox{if } u \in \cS_\rho, \\
				J'(u) & \quad \mbox{if } u \in \cD_\rho \setminus \cS_\rho.
			\end{array} \right.
		\end{align*}
		In particular, $J|_{\cD_\rho}'(u) = J'(u) + \lambda_u u$ for some $\lambda_u \in \R$ with $\lambda_u = 0$ if $u \in \cD_\rho \setminus \cS_\rho$.
	\end{Rem}
	
	To debate the dynamical properties related to \eqref{eq:time}, we need the following definition.
	\begin{Def}\label{def:os}
		We say that $X \subset H^1(\rn;\mathbb{C})$ is orbitally stable if and only if for every $\varepsilon > 0$ there exists $\delta > 0$ such that for any $\psi_0 \in H^1(\rn;\mathbb{C})$, if $\inf_{u \in X} \|u - \psi_0\| \le \delta$, then $\sup_{t > 0} \inf_{u \in X} \|u - \Psi(\cdot,t)\| \le \varepsilon$, where $\Psi$ is the unique solution to \eqref{eq:time}.
	\end{Def}
	\noindent We will mention the orbital stability only in situations where \eqref{eq:time} is \textit{globally} well-posed (cf. \cite[Definition 3.1.5]{Cazenave}), hence Definition \ref{def:os} makes sense.
	
	We have the following result that concerns all the local minimizers as in Theorem \ref{th:locmin}.
	\begin{Prop}\label{pr:os}
		If (\ref{F0})--(\ref{F2}) and \eqref{eq:rho} hold, $f$ is odd, and there exists $q \in (2,2^*)$ such that (at least) one of the following holds:
		\begin{enumerate}[label=(F5\alph{*}),ref=F5\alph{*}]
			\item \label{eq:os1} $|f(t) - f(s)| \lesssim (1 + |t| + |s|)^{q-2} |t - s|$ for all $ t,s \in \R$; 
			\item \label{eq:os2} $f \in \cC^1(\R)$ and $|f'(t)| \lesssim |t|^{q-2} + |t|^{2^*-2}$  for all $t \in \R$;
		\end{enumerate}
		then the solutions to \eqref{eq:time} exist globally in time and the set
		\begin{equation*}
			\cG := \Set{u \in \cU_{R_0}(\rho) | J(u) = m_{R_0}(\rho)}
		\end{equation*}
		is orbitally stable.
	\end{Prop}
	\noindent Notice that, if (\ref{eq:os2}) holds, then (\ref{F2}) implies that $2 < q < 2_\#$.
	
	When (\ref{eq:os1}) holds, Proposition \ref{pr:os} improves, to some extent, Soave's results \cite{Soave} with $N\ge3$ because, as it is clear from its proof, the threshold on $\rho$ to have a solution $(u,\lambda)$ to \eqref{eq:main} such that $J(u) = m_{R_0}(\rho)$ is the same as the one to obtain that $\cG$ is orbitally stable. In particular,
	for $q \in (2,2^*)$, denoting by $C_{N,q}$ the optimal constant in the Gagliardo--Nirenberg inequality, i.e.,
	\begin{equation}\label{eq:GN}
		|u|_q \le C_{N,q} |\nabla u|_2^{\gamma_q} |u|_2^{1-\gamma_q}, \quad u \in \hrn,
	\end{equation}
	where $\gamma_q = N(q-2)/(2q)$, we have the following outcome.

	\begin{Cor}
		Let $f(t) = |t|^{p-2}t + \mu |t|^{q-2}t$, $t \in \R$, with $\mu>0$ and $2 < q < 2_\# < p < 2^*$. If
		\begin{equation*}
			\left(\mu \rho^{q(1-\gamma_q)}\right)^{\gamma_pp-2} \rho^{p(1-\gamma_p)(2-\gamma_qq)} < \left(\frac{p(2-\gamma_qq)}{2C_{N,p}^p(\gamma_pp-\gamma_qq)}\right)^{2-\gamma_qq} \left(\frac{q(\gamma_pp-2)}{2C_{N,q}^q(\gamma_pp-\gamma_qq)}\right)^{\gamma_pp-2}
		\end{equation*}
		holds, then $\cG$ is orbitally stable.
	\end{Cor}

	Now we turn to the existence of a second solution to \eqref{eq:main}. Let us define
	\begin{equation}\label{eq:H}
		H(t) := f(t)t - 2F(t) \quad \text{for all } t \in \R.
	\end{equation}
	We assume that $H = H_1 + H_2$, where $H_1$ and $H_2$ satisfy what follows.

	\begin{enumerate}[label=(H\arabic{*}),ref=H\arabic{*}]\setcounter{enumi}{-1}
		\item \label{H0n} $H_1,H_2 \in \cC^1(\R;\R)$ and there exist $a \in (2, 2_\#)$ and $b \in (2_\#, 2^*)$ such that
		$$
		H_1(t) \lesssim |t|^2 + |t|^a, \quad H_2(t) \lesssim |t|^b + |t|^{2^*} \ \mbox{for all } t \in \R.
		$$
		\item \label{H2n} There holds 
		$$
		2 H_1(t) \leq h_1(t)t \leq a H_1(t), \quad b H_2(t) \leq h_2(t)t \leq 2^* H_2(t) \quad \text{for all } t \in \R,
		$$
		where $h_j := H_j'$ for $j \in \{1,2\}$.
	\end{enumerate}

	Observe that the inequalities in (\ref{H0n}), together with (\ref{H2n}), imply that $|h(t)| \lesssim |t|+|t|^{2^*-1}$ for all $t \in \R$, where $h=H'$.

	When $f$ is odd, we assume in addition that $H_j$ are even and, to deal with complex-valued functions, as for $F$, we define $H_j(z) := H_j(|z|)$ for all $z \in \mathbb{C}$ and $j \in \{1,2\}$, $H(z) := H_1(z)+H_2(z)$, and $h(z):=h_1(z)+h_2(z)$ for $z \in \mathbb{C}$. Since $H_j$ are even, $h_j$ are odd and we extend them to $\mathbb{C}$ as we do with $f$.
	In this case, (\ref{H2n}) implies
	\[
	2 H_1(z) \leq h_1(z) \overline{z} \leq a H_1(z), \quad b H_2(z) \leq h_2(z)\overline{z} \leq 2^* H_2(z) \quad \text{for all } z \in \mathbb{C}.
	\]
	
	In view of (\ref{F0}), every solution to \eqref{eq:main} satisfies the following Poho\v{z}aev identity (\cite{Poh})
	$$
	\int_{\R^N} |\nabla u|^2 \, \dx
	= 2^* \left(\int_{\R^N} F(u) \, \dx - \frac{\lambda}{2} \int_{\R^N}  |u|^2 \, \dx\right).
	$$
	Moreover, every solution satisfies the Nehari identity
	$$
	\int_{\R^N} |\nabla u|^2 \, \dx +\lambda \int_{\R^N} |u|^2 \, \dx - \int_{\R^N} f(u)u \, \dx = 0
	$$
	(where $f(u)u$ is replaced with $f(u)\overline{u}$ if $u$ is complex-valued). Taking a linear combination of these two identities, we find that every solution to \eqref{eq:main} satisfies
	\begin{equation}\label{defM}
		M(u):= \int_{\R^N} |\nabla u|^2 \, \dx - \frac{N}{2} \int_{\R^N} H(u) \, \dx = 0.
	\end{equation}
	
	We introduce following constraints
	\begin{equation*}
		\cM := \Set{u \in \hrn \setminus \{0\} | M(u)=0}, \quad \cM^\textup{rad} := \cM \cap H^1_\textup{rad}(\rn).
	\end{equation*}
	Note that
	\begin{equation}\label{equivM}
		u \in \cM \quad \text{ if and only if } \quad \left.\frac{\mathrm{d}}{\mathrm{d}s} J(s\star u)\right|_{s=1} = 0\text{ and }u \ne 0
	\end{equation}
	and that $\cM$ consists of the disjoint union of the following sets
	\begin{align*}
		\cM_0 & := \Set{u \in \cM | \left.\frac{\mathrm{d}^2}{\mathrm{d}s^2} J(s\star u)\right|_{s=1} = 0},\\
		\cM_- & := \Set{u \in \cM | \left.\frac{\mathrm{d}^2}{\mathrm{d}s^2} J(s\star u)\right|_{s=1} < 0},\\
		\cM_+ & := \cM \setminus \left(\cM_- \cup \cM_0\right),
	\end{align*}
	where, for $u\in\cM$,
	\begin{equation}
		\label{phisec}
		\left.\frac{\mathrm{d}^2}{\mathrm{d}s^2} J(s\star u)\right|_{s=1}
		=
		\frac{N^2}{4} \left( 2_\# \int_{\R^N} H(u) \, \dx -  \int_{\R^N} h(u)u \, \dx \right).
	\end{equation}
	We shall also consider the set $\cM_-^\textup{rad} := \cM_- \cap H^1_\textup{rad}(\rn)$.

	Our second existence result is the following.
	
	\begin{Th}\label{th:2sol}
		If (\ref{F1})--(\ref{F4}), (\ref{H0n}), (\ref{H2n}), and (\ref{J1}) hold and $\rho>0$ is sufficiently small, then there exist $\widetilde{u} \in \cS_\rho$ and $\lambda_{\widetilde u} > 0$ such that $J(\widetilde{u}) = \min_{\cM_-^\textup{rad} \cap \cD_\rho} J > 0$ and $(\widetilde{u},\lambda_{\widetilde u})$ is a solution to \eqref{eq:main}. Moreover, if $f$ is odd or $f|_{(-\infty,0)} \equiv 0$, then $J(\widetilde{u}) = \min_{\cM_- \cap \cD_\rho} J$ and $\widetilde u$ can be chosen to be positive and non-increasing in the radial coordinate; in this case, $\widetilde u$ is actually decreasing if there exists $\tau > 0$ such that $f(t) \le \lambda_{\widetilde u} t$ for all $t \in [0,\tau]$ and $f(t) > \lambda_{\widetilde u} t$ for all $t > \tau$.
	\end{Th}
	
	The next definition and proposition are somewhat antipodal to the stability results mentioned above.
	
	\begin{Def}
		We say that $u \in H^1(\rn;\mathbb{C})$ is strongly unstable if and only if for every $\varepsilon > 0$ there exists $\psi_0 \in H^1(\rn; \mathbb{C})$ such that $\|u-\psi_0\| < \varepsilon$ and $|\nabla \Psi|_2$ blows up in finite time, where $\Psi$ is the unique solution to \eqref{eq:time}.
	\end{Def}
	
	\begin{Prop}\label{pr:si}
		If (\ref{J2}), (\ref{eq:os1}), and the assumptions of Theorem \ref{th:2sol} hold and $f$ is odd, then $\widetilde{u}$ is strongly unstable.
	\end{Prop}

	We conclude the Introduction with a list of open problems.

	\begin{itemize}
		\item[(OP1)]
		In \cite[Theorem 1.4]{JL_min}, see also \cite{Maris}, in a context where $J|_{\cS_\rho}$ is bounded below, it was proved that every global minimizer is radial (about a point) and radially non-increasing without assuming that $f$ is odd. In the context of this paper, where global minimizers do not exist, are all the \textit{local minimizers} as in Theorem \ref{th:locmin} radial and radially non-increasing? Or -- at least -- is this true if they are also least-energy solutions (cf. Proposition \ref{pr:gs})?
		
		\item[(OP2)]
		Assumption (\ref{eq:os1}) implies that the growth of $f$ at infinity is controlled by a Sobolev-subcritical power. Can this assumption be replaced with simply the Sobolev-subcritical growth? In other words, can (\ref{eq:os1}) be replaced with
		\[
		|f(t) - f(s)| \le \mathcal{L}(r) |t-s| \quad \text{for all } t,s \in [-r,r] \text{ and all } r > 0,
		\]
		where $\mathcal{L} \colon (0,+\infty) \to (0,+\infty)$ is such that $\lim_{r \to +\infty} \mathcal{L}(r) / r^{2^*-2} = 0$?
		
		\item[(OP3)]
		Assumption (\ref{eq:os2}) implies that the growth of $f'$ at the origin is controlled by a positive power. Can this assumption be replaced with simply the convergence to $0$? In other words, can the second part of (\ref{eq:os2}) be replaced with
		\[
		|f'(t)| \lesssim \mathcal{L}(t) + |t|^{2^*-2} \quad \text{for all } t \in \R,
		\]
		where $\mathcal{L} \colon (0,+\infty) \to (0,+\infty)$ is such that $\lim_{t \to 0^+} \mathcal{L}(t) = 0$?
	\end{itemize}
	
	In what follows, $C$ denotes a generic positive constant which may vary from one line to another. When $f$ is odd, it is understood that, for $v \in H^1(\R^N;\mathbb{C})$, $f(v)v$ and $h(v)v$ are replaced with $f(v)\overline{v}$ and $h(v)\overline{v}$ respectively.

	\section{Local minimum and the first normalized solution}\label{sec:lm}
	
	The aim of the section is to study the existence of negative-energy solutions (Theorem \ref{th:locmin} and Proposition \ref{pr:gs}) and their orbital stability (Proposition \ref{pr:os}).
	
	\subsection{The elliptic problem}
	
	Let us start with the proof of the existence of negative-energy solutions by finding a proper estimate below of the energy functional $J$ on the disk $\cD_\rho$. 
	
	Observe that, by \eqref{eq:C0} and \eqref{eq:gns}, for $u \in \cD_\rho$,
	\begin{equation}\label{eq:Jbddb}
		J(u) \ge \frac12 |\nabla u|_2^2 - C_0 \left(|u|_2^2 + |u|_{2^*}^{2^*}\right) \ge \frac12 |\nabla u|_2^2 - C_0 \rho^2 - C_0 S^{-2^* / 2} |\nabla u|_2^{2^*}
	\end{equation}
	and so $J|_{\cD_\rho}$ is bounded below over bounded subsets in $\hrn$ or, equivalently, $D^{1,2}(\rn)$.\\
	With this in mind, define $g \colon (0,+\infty) \times (0,+\infty) \to \R$ as
	\begin{equation}\label{eq:g}
		g(\alpha,t) := \frac{1}{2} - C_0 \alpha^2 t^{-2} - C_0 S^{-2^* / 2} t^{2^*-2}.
	\end{equation}
	Then, \eqref{eq:Jbddb} reads as
	\begin{equation}\label{eq:Jg}
		J(u) \ge g(\rho,|\nabla u|_2) |\nabla u|_2^2 \quad
		\text{for all }
		u \in \cD_\rho
	\end{equation}
	\begin{Lem}\label{le:g}
		The following facts hold.
		\begin{enumerate}[label=(g$_\arabic{*}$),ref=g$_\arabic{*}$]
			\item \label{max} For every $\alpha>0$, the function $t \mapsto g(\alpha,t) t^2$ has a unique critical point, which is a global maximizer.
			\item \label{eq:g_prop} If \eqref{eq:rho} holds, then there exist $R_0, R_1>0$, $R_0 < R_1$, such that $g(\rho, R_0) = g(\rho, R_1) = 0$, $g(\rho, t) > 0$ for $t \in (R_0,R_1)$, and $g(\rho, t) < 0$ for $t \in (0,R_0) \cup (R_1, +\infty)$.
			\item \label{le:g-ineq} If $t>0$ and $\alpha_1 \ge \alpha_2 > 0$, then for every $s \in \left[t\alpha_2/\alpha_1,t\right]$ there holds $g(\alpha_2,s) \ge g(\alpha_1,t)$.
			\item \label{g4} If $\rho$ satisfies \eqref{eq:rho}, then there exists $\varepsilon > 0$ such that \eqref{eq:rho} is verified by every $\rho' \in (\rho-\varepsilon, \rho+\varepsilon)$ and the functions $(\rho - \varepsilon, \rho + \varepsilon) \ni \rho' \mapsto R_i(\rho ') \in (0,+\infty)$, $i \in \{0,1\}$ with $R_i$ defined in (\ref{eq:g_prop}), are invertible and of class $\cC^1$.
		\end{enumerate}
	\end{Lem}
	\begin{proof}
		Points (\ref{max}) and (\ref{eq:g_prop}) follow from direct computations.\\
		Concerning (\ref{le:g-ineq}), it is clear that $g(\alpha_2, t) \ge g(\alpha_1, t)$ for all $t>0$. Moreover,
		\[
		g\left(\alpha_2, \frac{\alpha_2}{\alpha_1}t\right) - g(\alpha_1, t) = \frac{C_0}{S^{2^* / 2}} \left[1-\left(\frac{\alpha_2}{\alpha_1}\right)^{2^*-2} \right] t^{2^*-2} \ge 0,
		\]
		hence we conclude thanks to (\ref{max}) and the regularity of $g(\alpha_2,\cdot)$.\\
		Finally, as for (\ref{g4}), it follows from the regularity of $g$, $\partial_t g(\rho, R_0) > 0$, $\partial_t g(\rho, R_1) < 0$, and the implicit function theorem.
	\end{proof}

	Now, using that, for $u \in \hrn \setminus \{0\}$ and $s>0$, the operator $\star$ defined in \eqref{star} satisfies
	\[
	|s \star u|_2 = |u|_2 \quad
	\text{and} \quad
	|\nabla (s \star u)|_2 = s |\nabla u|_2.
	\]
	Let us show that, in our setting of assumptions, the infimum of $J$ in $\mathcal{U}_r(\alpha)$  (see \eqref{URamRa}) is negative.
	
	\begin{Lem}\label{le:neg}
		If (\ref{F0}) and (\ref{F2}) hold, then $m_R(\alpha) \in (-\infty,0)$ for every $\alpha, R > 0$.
	\end{Lem}
	\begin{proof}
		First of all, let us observe that, from \eqref{eq:Jbddb}, $m_R(\alpha) > -\infty$. Then, let us fix $u \in \cD_\alpha \cap L^\infty(\rn) \setminus \{0\}$ and observe that, for $s > 0$,
		\begin{equation*}
			J(s \star u) = s^2 \left(\frac{|\nabla u|_2^2}{2} - \frac{1}{\left(s^{N/2}\right)^{2_\#}}\int_{\rn} F\big(s^{N/2}u\big) \, \dx\right).
		\end{equation*}
		Note that, from (\ref{F2}), $F(s^{N/2}u) > 0$ a.e. in $\mathrm{supp} \, u$ for sufficiently small $s > 0$. Therefore, from Fatou's Lemma and again (\ref{F2}),
		\[
		\lim_{s \to 0^+} \frac{1}{\left(s^{N/2}\right)^{2_\#}}\int_{\rn} F\big(s^{N/2}u\big) \, \dx = +\infty.
		\]
		This implies that $J(s \star u) < 0$ for sufficiently small $s > 0$, and since $s \star u \in \cU_R(\alpha)$ provided $s > 0$ is small, we can conclude. 
	\end{proof}

	Now we show a property ensuring that $m_{R_0}(\rho)$ cannot be attained at the boundary of $\mathcal{U}_{R_0}(\rho)$.
	\begin{Rem}\label{re:bdda}
		From (\ref{eq:g_prop}) and Lemma \ref{le:neg}, there exists $\varepsilon > 0$ such that
		\[
		0 \ge g(\rho,s) \ge \frac{m_{R_0}(\rho)}{2 R_0^2}
		\text{ for all } s \in [R_0-\varepsilon,R_0].
		\]
		This, \eqref{eq:Jg}, and Lemma \ref{le:neg} yield that for all $u \in \cD_\rho$ such that $R_0 - \varepsilon \le |\nabla u|_2 \le R_0$ there holds
		\[
		J(u)
		\ge g(\rho,|\nabla u|_2) |\nabla u|_2^2 
		\ge R_0^{2} \frac{m_{R_0}(\rho)}{2 R_0^2} > m_{R_0}(\rho).
		\]
	\end{Rem}
	
	Next, we show the following \textit{subadditivity} property.
	\begin{Lem}\label{le:sub}
		If (\ref{F0}), (\ref{F2}), and \eqref{eq:rho} hold, then for every $\alpha \in (0,\rho)$ there holds
		\[
		m_{R_0}(\rho) \le m_{R_0}(\alpha) + m_{R_0}\left(\sqrt{\rho^2-\alpha^2}\right).
		\]
		Moreover, the inequality above is strict if at least one between $m_{R_0}(\alpha)$ and $m_{R_0}\left(\sqrt{\rho^2-\alpha^2}\right)$ is attained.
	\end{Lem}
	\begin{proof}
		We first prove that
		\begin{equation}\label{eq:theta}
			m_{R_0}(\theta \alpha) \le \theta^2 m_{R_0}(\alpha)
		\end{equation}
		for every $\theta \in \left[1,\rho/\alpha\right]$ and that the inequality is strict if $m_{R_0}(\alpha)$ is attained and $\theta>1$.\\
		From \eqref{eq:Jg}, Lemma \ref{le:neg}, and the definition of $m_{R_0}(\alpha)$, for every $\varepsilon > 0$ sufficiently small there exists $u \in \cU_{R_0}(\alpha)$ such that
		\begin{equation}\label{gal0}
			g(\alpha,|\nabla u|_2) |\nabla u|_2^2 \le J(u) \le m_{R_0}(\alpha) + \varepsilon < 0.
		\end{equation}
		If $|\nabla u|_2 \ge R_0 \alpha/\rho$, then from (\ref{le:g-ineq}) with $\alpha_1 = \rho$, $\alpha_2 = \alpha$, $t = R_0$, and $s = |\nabla u|_2$ we obtain
		\[
		g(\alpha,|\nabla u|_2) \ge g(\rho,R_0)=0,
		\]
		which is in contradiction with \eqref{gal0}. Therefore $u \in \cU_{\alpha R_0/\rho}(\alpha)$. Define $v := u\left(\cdot/\theta^{2/N}\right) \in \cU_{R_0}(\theta \alpha)$. Then
		\begin{equation*}
			m_{R_0}(\theta \alpha) \le J(v) = \theta^2 \left(\frac{1}{2 \theta^{4/N}}|\nabla u|_2^2 - \int_{\rn} F(u) \, \dx\right) \le \theta^2 J(u) \le \theta^2 \left( m_{R_0}(\alpha) + \varepsilon \right)
		\end{equation*}
		and, by the arbitrariness of $\varepsilon$, we get \eqref{eq:theta}.\\
		In addition, if $m_{R_0}(\alpha)$ is attained and $\theta>1$, we can take $u\in \cU_{R_0}(\alpha)$ as the minimizer and, repeating the previous argument, we have
		\begin{equation}\label{eq:thetaSTR}
			m_{R_0}(\theta \alpha) < \theta^2 J(u) = \theta^2 m_{R_0}(\alpha).   
		\end{equation}
		If $\alpha^2 > \rho^2 - \alpha^2$, we obtain
		\begin{equation*}
			m_{R_0}(\rho)
			\overset{(1)}{\le} \frac{\rho^2}{\alpha^2} m_{R_0}(\alpha) 
			= m_{R_0}(\alpha) + \frac{\rho^2 - \alpha^2}{\alpha^2} m_{R_0}(\alpha)
			\overset{(2)}{\le} m_{R_0}(\alpha) + m_{R_0}\left(\sqrt{\rho^2 - \alpha^2}\right),
		\end{equation*}
		where the inequality $(1)$ (respectively, $(2)$) is strict if $m_{R_0}(\alpha)$ (respectively, $m_{R_0}\left(\sqrt{\rho^2 - \alpha^2}\right)$) is attained.\\
		Analogously, if $\alpha^2 < \rho^2 - \alpha^2$, we obtain
		\begin{align*}
			m_{R_0}(\rho) 
			& \overset{(3)}{\le}
			\frac{\rho^2}{\rho^2 - \alpha^2} m_{R_0}\left(\sqrt{\rho^2 - \alpha^2}\right)
			= m_{R_0}\left(\sqrt{\rho^2 - \alpha^2}\right) + \frac{\alpha^2}{\rho^2 - \alpha^2} m_{R_0}\left(\sqrt{\rho^2 - \alpha^2}\right)\\
			& \overset{(4)}{\le} m_{R_0}\left(\sqrt{\rho^2 - \alpha^2}\right) + m_{R_0}(\alpha),
		\end{align*}
		where the inequality $(3)$ (respectively, $(4)$) is strict if $m_{R_0}\left(\sqrt{\rho^2 - \alpha^2}\right)$ (respectively, $m_{R_0}(\alpha)$) is attained.\\
		Finally, if $\alpha^2 = \rho^2 - \alpha^2$ (i.e. $\sqrt{2} \, \alpha = \rho$), then
		\[
		m_{R_0}( \sqrt{2} \, \alpha ) \le 2m_{R_0}(\alpha)
		\]
		follows from \eqref{eq:theta}, while the inequality being strict
		provided $m_{R_0}(\alpha)$ is attained follows from \eqref{eq:thetaSTR}.
	\end{proof}
	
	The next Lemma shows the relative compactness in $L^p(\R^N)$ of minimizing sequences.
	
	\begin{Lem}\label{le:minim}
		Assume that (\ref{F0})--(\ref{F2}) and \eqref{eq:rho} hold. If $\rho_n \to \rho$, $R_n \to R_0$, and
		$(u_n) \subset H^1(\R^N)$ is such that $u_n\in \cU_{R_n}(\rho_n)$ for every $n\in\mathbb{N}$ and $J(u_n) \to m_{R_0}(\rho)$, then there exist $(y_n) \subset \rn$ and $\overline{u} \in \cD_\rho \setminus \{0\}$ such that $|\nabla \overline{u}|_2 \le R_0$ and $u_n(\cdot + y_n) \to \overline{u}$ in $L^p(\rn)$ for every $p \in [2,2^*)$.
	\end{Lem}
	
	\begin{proof}
		First of all, since $(u_n)$ is bounded in $H^1(\R^N)$, it suffices to prove the statement for $p=2$ in virtue of the interpolation inequality and the Sobolev embedding.\\
		Assume by contradiction that
		\begin{equation*}
			\lim_{n} \sup_{y\in\rn} \int_{B(y,1)} |u_n|^2 \, \dx = 0.
		\end{equation*}
		Then, from Lions' Lemma \cite[Lemma I.1]{Lions} we obtain that $\lim_n |u_n|_q = 0$ for every $q \in (2,2^*)$.\\
		Fix $q \in (2,2^*)$ and $\varepsilon > 0$. From (\ref{F1}) and \eqref{eq:C0}, there exists $c = c(q,\varepsilon,C_0) > 0$ such that for every $t \in \R$
		\[
		F(t)\leq \varepsilon t^2 + c |t|^q + C_0 |t|^{2^*}.
		\]
		Recalling that $g(\rho,R_0) = 0$, for $\varepsilon \ll 1$, using also Lemma \ref{le:neg}, there holds
		\begin{align*}
			0
			&
			> m_{R_0}(\rho)
			= \lim_n J(u_n)
			\ge \limsup_n \left(\frac{1}{2}|\nabla u_n |_2^2 - \varepsilon |u_n|_2^2 - c |u_n|_q^q - C_0 |u_n|_{2^*}^{2^*}\right)\\
			&
			\ge \left(\frac12 - \frac{C_0 R_0^{2^*-2}}{S^{2^*/2}}  \right) \limsup_n |\nabla u_n|_2^2 
			- \varepsilon \rho^2 
			= \frac{C_0 \rho^2}{R_0^2}\limsup_n |\nabla u_n|_2^2 - \varepsilon \rho^2
			\geq - \varepsilon \rho^2,
		\end{align*}
		a contradiction.
		Hence there exists $(y_n) \subset \R^N$ such that
		$$
		\limsup_n \int_{B(0,1)} |u_n(\cdot + y_n)|^2 \, \dx > 0.
		$$
		As a consequence, there exists $\overline u \in \cD_\rho\setminus\{0\}$ with $|\nabla \overline{u}|_2\leq R_0$
		such that, setting $v_n(x) := u_n(x+y_n) - \overline u(x)$, up to a subsequence $v_n \rightharpoonup 0$ in $\hrn$ and $v_n(x) \to 0$ for a.e. $x \in \rn$.\\
		Denote $\alpha := |\overline{u}|_2$ and assume by contradiction that $\beta := \lim_n |v_n|_2 > 0$ along a subsequence. From (\ref{F0}), for every $\delta>0$ there exists $C_\delta>0$ such that for all $t,s \in \R$ there holds
		\[
		|F(t+s) - F(t)| \le \delta (t^2 + |t|^{2^*}) + C_\delta (s^2 + |s|^{2^*}).
		\]
		To see this, we compute
		\begin{align*}
			|F(t+s) - F(t)| & \le \int_0^1 |f(t + \tau s)| |s| \, \mathrm{d}\tau \le C \int_0^1 \left(|t+\tau s| + |t + \tau s|^{2^*-1}\right) |s| \, \mathrm{d}\tau\\
			& \le C \left(|t| + |s| + |t|^{2^*-1} + |s|^{2^*-1}\right) |s| \le \delta (t^2 + |t|^{2^*}) + C_\delta (s^2 + |s|^{2^*}),
		\end{align*}
		where the last step comes from Young's inequality. Then, from the Brezis--Lieb Lemma \cite[Theorem 2]{BL} and the weak convergence in $\hrn$ we get
		\[
		\lim_n \big[J(u_n) - J(v_n)\big] = J(\overline u),
		\quad
		\lim_n \big[|u_n|_2^2 - |v_n|_2^2 \big]= |\overline u|_2^2,
		\quad
		\lim_n \big[|\nabla u_n|_2^2 - |\nabla v_n|_2^2 \big]= |\nabla \overline u|_2^2,
		\]
		and so, in particular, $\beta^2\leq \rho^2-\alpha^2$ and $|\nabla v_n|_2 < R_0$ for $n \gg 1$.\\
		Now, for every $n$ define
		\[
		\widetilde{v}_n :=
		\begin{cases}
			v_n & \text{if } |v_n|_2 \le \beta,\\
			\beta v_n /|v_n|_2& \text{if } |v_n|_2 > \beta.
		\end{cases}
		\]
		Then $\widetilde{v}_n \in \cU_{R_0}(\beta)$ and, arguing as in \cite[Lemma 2.4]{Shibata14} and using (\ref{F0}), we obtain that
		$$\lim_n \big[J(v_n) - J(\widetilde{v}_n) \big]= 0.$$
		Hence, from Lemma \ref{le:sub}, the continuity of $J$, and since $m_{R_0}$ is non-increasing, up to a subsequence there holds
		\begin{align*}
			m_{R_0}(\rho) & = \lim_n J(u_n) = J(\overline u) + \lim_n J(\widetilde{v}_n) \ge m_{R_0}(\alpha) + m_{R_0}(\beta)\\
			& \ge m_{R_0}\left(\sqrt{\alpha^2+\beta^2}\right) \ge m_{R_0}(\rho).
		\end{align*}
		This implies that $m_{R_0}(\alpha)$ is attained at $\overline u$ and so, again from Lemma \ref{le:sub},
		\[
		m_{R_0}(\alpha) + m_{R_0}(\beta) > m_{R_0}\left(\sqrt{\alpha^2+\beta^2}\right),
		\]
		a contradiction.
	\end{proof}
	
	We are ready to prove the first main result.
	
	\begin{proof}[Proof of Theorem \ref{th:locmin}]
		Let $(\rho_n) \subset (0,+\infty)$, $(R_n) \subset (0,+\infty)$, and $(u_n) \subset H^1(\R^N)$ be as in the statement.
		From Lemma \ref{le:minim}, there exists $\overline{u} \in \cD_{\rho} \setminus \{0\}$ with $|\nabla \overline{u}|_2 \le R_0$ such that, replacing $u_n(\cdot + y_n)$ with $u_n$ and up to a subsequence,
		$u_n \rightharpoonup \overline{u}$ in $\hrn$ and $u_n \to \overline{u}$ in $L^2(\rn)$ and a.e. in $\R^N$, and, by Remark \ref{re:bdda},
		$J(\overline u) \ge m_{R_0}(\rho)$.\\
		Denote $v_n := u_n - \overline u$. From the weak convergence,
		$$\lim_n \big(|\nabla u_n|_2^2 - |\nabla v_n|_2^2\big) = |\nabla \overline{u}|_2^2 > 0,$$
		thus $|\nabla v_n|_2 < R_0$ for $n \gg 1$. Moreover, from the Bresiz--Lieb Lemma \cite[Theorem 2]{BL},
		\[
		m_{R_0}(\rho) = \lim_n J(u_n) = \lim_n J(v_n) + J(\overline u) \ge \lim_n J(v_n) + m_{R_0}(\rho),
		\]
		which, together with $g(\rho, R_0) = 0$, \eqref{eq:C0}, and $\lim_n |v_n|_2 = 0$, implies
		\begin{align*}
			0
			& \ge
			\lim_n J(v_n) \ge \lim_n \left(\frac12 |\nabla v_n|_2^2 - C_0 \left(|v_n|_2^2 + |v_n|_{2^*}^{2^*}\right)\right)\\
			& \ge
			\left(\frac12 - \frac{C_0  R_0^{2^*-2}}{S^{2^* / 2}}\right) \lim_n |\nabla v_n|_2^2
			= \frac{C_0 \rho^2}{R_0^2} \lim_n |\nabla v_n|_2^2
			\ge 0,
		\end{align*}
		i.e., $\lim_n |\nabla v_n|_2 = 0$, whence $u_n \to \overline u$ in $\hrn$. It follows then that $J(\overline u) = \lim_n J(u_n) = m_{R_0}(\rho) < 0$
		(from Lemma \ref{le:neg}) and so, from Remark \ref{re:bdda}, $|\nabla \overline u|_2 < R_0$. 
		Consequently\footnote{The argument from this point on works for every $u$ as in the statement.}, there exists $\lambda_{\overline{u}} \in \R$ such that
		\[
		-\Delta \overline u + \lambda_{\overline{u}} \overline{u} = f(\overline u)
		\quad
		\text{in }\mathbb{R}^N.
		\]
		If $\lambda_{\overline{u}} \le 0$, then from the Poho\v{z}aev identity we obtain
		\[
		(N-2) \int_{\rn} |\nabla \overline u|^2 \, \dx = N \int_{\rn} \left(2 F(\overline u) - \lambda_{\overline{u}} |\overline u|^2\right) \, \dx \ge 2N \int_{\rn} F(\overline u) \, \dx
		\]
		and so $0 > J(\overline{u}) \ge |\nabla \overline{u}|_2^2/N$, which is a contradiction. This implies that $\lambda_{\overline{u}} > 0$ and, consequently, $\overline{u} \in \cS_{\rho}$ (otherwise we would have $\lambda_{\overline{u}} = 0$).\\
		Now we prove that $\overline u$ has constant sign.\\
		Let $\overline{u}_\pm := \max\{\pm\overline{u},0\}$ and $\rho_\pm := |\overline{u}_\pm|_2$. Of course, $\pm \overline{u}_\pm \in \cU_{R_0}(\rho_\pm) \subset \cU_{R_0}(\rho)$. If $|\overline{u}_+|_2 |\overline{u}_-|_2 > 0$, then, from Lemma \ref{le:sub} and since $\rho^2=\rho_+^2 + \rho_-^2$,
		\begin{equation*}
			m_{R_0}(\rho)
			= J(\overline u)
			= J(\overline{u}_+) + J(-\overline{u}_-)
			\ge m_{R_0}(\rho_+) + m_{R_0}(\rho_-)
			\ge m_{R_0}(\rho).
		\end{equation*}
		Hence $m_{R_0}(\rho_+)$ is attained at $\overline{u}_+$, $m_{R_0}(\rho_-)$ is attained at $-\overline{u}_-$, and, again from Lemma \ref{le:sub}, $m_{R_0}(\rho_+) + m_{R_0}(\rho_-) > m_{R_0}(\rho)$, a contradiction.\\
		Finally, writing
		\[
		-\Delta \overline u + \left(\lambda_{\overline u} + \max\left\{-\frac{f(\overline u)}{\overline u},0\right\}\right) \overline u = \max\left\{\frac{f(\overline u)}{\overline u},0\right\} \overline u,
		\]
		the fact that $\overline u$ does not vanish follows from the maximum principle \cite[Theorem 8.19]{GilTru} due to Remark \ref{rem:12} (\ref{Linfty}).
	\end{proof}

	\begin{proof}[Proof of Proposition \ref{pr:gs}]
		We shall prove that $m_{R_0}(\rho) = \inf \Set{J(u) | u \in \cD_\rho \text{ and } J|_{\cD_\rho}'(u) = 0}$, then Theorem \ref{th:locmin} will yield that such an infimum is actually a minimum.\\
		Since, from Theorem \ref{th:locmin}, $m_{R_0}(\rho)$ is attained, we have $m_{R_0}(\rho) \ge \inf \Set{J(u) | u \in \cD_\rho \text{ and } J|_{\cD_\rho}'(u) = 0}$.\\
		Assume by contradiction that the strict inequality holds, so there exists $u \in \cD_\rho \setminus \{0\}$ such that $J|_{\cD_\rho}'(u) = 0$ and $J(u) < m_{R_0}(\rho)$.
		From the definition of $m_{R_0}(\rho)$, there holds $|\nabla u|_2 \ge R_0$.
		In fact, since $J(u) < 0$, we know from \eqref{eq:Jg} and (\ref{eq:g_prop}) that $|\nabla u|_2 > R_1$.\\
		Consider the function $\varphi\colon(0,+\infty) \to \R$, $\varphi(s):=J(s\star u)$. Let us recall that, since every critical point of $J|_{\cD_\rho}$ belongs to $\cM$, $\varphi'(1)=0$ by \eqref{equivM}.  Again from \eqref{eq:Jg} and (\ref{eq:g_prop}), $\varphi$ is positive on $(R_0/|\nabla u|_2,R_1/|\nabla u|_2)$ and negative at $s=1$.
		Moreover, using (\ref{Linfty}) in Remark \ref{rem:12} and arguing as in the proof of Lemma \ref{le:neg}, $\varphi(s)<0$ for $s \ll 1$. Thus $\varphi$ has a local maximum point $t_u \in (0,1)$.
		From (\ref{J1}), $\varphi' < 0$ in a right-hand neighbourhood of $t_u$, hence, from (\ref{J2}), $\varphi' < 0$ in $(t_u,+\infty)$, in contradiction with $\varphi'(1) = 0$.
	\end{proof}
	
	\begin{Rem}
		If (\ref{F0})--(\ref{F2}), (\ref{J1}), and \eqref{eq:rho} hold, then the thesis of Proposition \ref{pr:gs} is still true provided $F\ge0$, $\lim_{|t| \to +\infty} F(t) / |t|^{2_\#} = +\infty$, and for every $u \in \cD_\rho$ the function $(t_u,+\infty) \ni s \mapsto J(s\star u) \in \R$ does not have saddle points, where $t_u>0$ is defined in (\ref{J1}). As a matter of fact, in view of Fatou's Lemma, $\varphi(s)\to-\infty$ as $s\to +\infty$, so $1$ can be neither a saddle point of $\varphi$, nor a maximum point -- due to (\ref{J1}), nor a minimum point (otherwise $\varphi$ should have an additional local maximum point), which is a contradiction. We point out that a sufficient condition for the function $(t_u,+\infty) \ni s \to J(s\star u) \in \R$ not to have saddle points is that $\cM_0 \cap \cD_\rho = \emptyset$, cf. Section \ref{sec:pes}.
	\end{Rem}

	\subsection{Dynamics}
	
	For the reader's convenience, we begin this subsection with some standard notations and definitions.

	For $p,r \in [1,+\infty]$ and $I \subset \R$ an interval, we define the space
	\begin{equation*}
		Y(p,r,I) := L^p (I; L^r (\R^N; \mathbb{C})),
	\end{equation*}
	where the norm is
	\[
	\|\Psi\|_{Y(p,r,I)} :=
	\begin{cases}
		\displaystyle \left( \int_I |\Psi(\cdot,t)|_r^p \, \dt \right)^{1/p} & \text{if } p < +\infty,\vspace{3pt}\\
		\esssup_{t \in I} |\Psi(\cdot,t)|_r & \text{if } p = +\infty.
	\end{cases}
	\]
	Similarly, we define
	\begin{equation*}
		X(p,r,I) := L^p (I; W^{1,r} (\R^N; \mathbb{C})),
	\end{equation*}
	where the norm is
	\[
	\|\Psi\|_{X(p,r,I)}
	:=
	\begin{cases}
		\displaystyle \left( \int_I \|\Psi(\cdot,t)\|_{W^{1,r}}^p \, \dt \right)^{1/p} & \text{if } p < +\infty,\vspace{3pt}\\
		\esssup_{t \in I} \|\Psi(\cdot,t)\|_{W^{1,r}} & \text{if } p = +\infty,
	\end{cases}
	\]
	and $\| u \|_{W^{1,r}} := \left( |\nabla u|_r^r + |u|_r^r \right)^{1/r}$ is the usual norm in $W^{1,r}(\R^N; \mathbb{C})$. If $I=(0,T)$ for some $T>0$, we denote
	$$
	X(p,r,T) := X(p,r,(0,T)), \quad Y(p,r,T) := Y(p,r,(0,T)).
	$$
	
	Suppose that (\ref{F0}) holds and $0 \in I$. We say that $\Psi$ is a \textit{weak solution} to \eqref{eq:time} on $I$(cf. \cite[Definition 3.1.1]{Cazenave}) if and only if
	\begin{itemize}
		\item $\Psi \in X(\infty, 2, I) \cap W^{1,\infty} (I; H^{-1}(\R^N; \mathbb{C}))$,
		\item $\Psi$ satisfies \eqref{eq:time} in $H^{-1}(\R^N; \mathbb{C})$ for a.e. $t \in I$,
		\item $\Psi(0) = \psi_0$.
	\end{itemize}

	If $\Psi \in \cC(I; H^1(\R^N; \mathbb{C})) \cap \cC^1 (I; H^{-1} (\R^N; \mathbb{C}))$ satisfies \eqref{eq:time} in $H^{-1} (\R^N; \mathbb{C})$ for all $t \in I$ and $\Psi(0) = \psi_0$, then $\Psi$ is called a \textit{strong solution} to \eqref{eq:time} on $I$.

	We recall from \cite[Definition 3.1.5]{Cazenave} that the time-dependent problem \eqref{eq:time} is \textit{locally well-posed} in $H^1(\R^N; \mathbb{C})$ if and only if:
	\begin{enumerate}[label=(\roman{*}),ref=\roman{*}]
		\item \textit{the uniqueness property} in $H^1(\R^N; \mathbb{C})$ of \eqref{eq:time} holds, i.e., for any initial condition $\psi_0 \in H^1(\R^N; \mathbb{C})$ and any time interval $I$ with $0 \in I$, any two weak solutions to \eqref{eq:time} on $I$ coincide;
		\item for every initial condition $\psi_0 \in H^1(\R^N; \mathbb{C})$ there is a strong solution $\Psi$ defined on a maximal interval of existence $(T_{\min} (\psi_0), T_{\max} (\psi_0)) {\color{blue}\ni 0}$;
		\item \textit{the blow-up alternative holds}: if $T_{\max} (\psi_0) < +\infty$, then $\lim_{t \to T_{\max}^{-} (\psi_0)} \| \Psi(t) \| = +\infty$ (and similarly for $T_{\min}(\psi_0)$);
		\item if $\psi_n \to \psi_0$ in $H^1 (\R^N; \mathbb{C})$ and $I \subset (T_{\min} (\psi_0), T_{\max} (\psi_0))$ is a closed interval, then the maximal solution $\Psi_n$ of \eqref{eq:time} with an initial condition $\Psi_n(0) = \psi_n$ is defined on $I$ for sufficiently large $n$, and $\Psi_n \to \Psi_0$ in $\cC (I; H^1(\R^N; \mathbb{C}))$.
	\end{enumerate}
	
	Here and in what follows, the gradient $\nabla$ is always meant with respect to the variable $x \in \rn$. Observe that both $Y(p,r,T)$ and $X(p,r,T)$ are Banach spaces (see, e.g., \cite[p. 97]{DiestelUhl}).

	We now prove Proposition \ref{pr:os}, starting with the Sobolev-subcritical case.
	
	\begin{proof}[Proof of Proposition \ref{pr:os} if (\ref{eq:os1}) holds]
		From \cite[Corollary 4.3.3]{Cazenave}, the Cauchy Problem \eqref{eq:time} is locally well-posed and there is conservation of mass and energy. Assume by contradiction that there exist $(\psi_n)_n \subset H^1(\R^N;\mathbb{C})$ and $(t_n)_n \subset (0,+\infty)$ such that
		\begin{equation*}
			\lim_n \inf_{u \in \cG} \|u - \psi_n\| = 0 \quad \text{and} \quad \limsup_n \inf_{u \in \cG} \|u - \Psi_n(\cdot,t_n)\| > 0,
		\end{equation*}
		where $\Psi_n$ is the unique solution to \eqref{eq:time} with $\psi_0 = \psi_n$. In particular, since $m_{R_0}(\rho)$ is attained at some $\overline u \in \cU_{R_0}(\rho) \cap \cS_\rho$ in virtue of Theorem \ref{th:locmin}, we have $\lim_n |\psi_n|_2 = \rho$, $\lim_n |\nabla \psi_n|_2 = |\nabla \overline u|_2 < R_0$, and $\lim_n J(\psi_n) = m_{R_0}(\rho) < 0$. From the conservation of mass and energy, for every $t>0$ such that $\Psi_n(\cdot,t)$ exists (in particular, for $t=t_n$), there holds $\lim_n |\Psi_n(\cdot,t)|_2 = \rho$ and $\lim_n J\bigl(\Psi_n(\cdot,t)\bigr) = m_{R_0}(\rho)$. In particular, if there exists $s_n>0$ such that $|\nabla \Psi_n(\cdot,s_n)|_2 \ge R_{0,n}$, where $R_{0,n} > 0$ is given by (\ref{eq:g_prop}) with $|\psi_n|_2$ instead of $\rho$, then there exists $r_n \in (0,s_n]$ such that $|\nabla \Psi_n(\cdot,r_n)|_2 = R_{0,n}$. Whence, from \eqref{eq:Jg}, $J\bigl(\Psi_n(\cdot,r_n)\bigr) \ge |\nabla \Psi_n(\cdot,r_n)|_2^2 g(|\Psi_n(\cdot,r_n)|_2, |\nabla \Psi_n(\cdot,r_n)|_2) = 0$ (observe that $|\Psi_n(\cdot,r_n)|_2$ satisfies \eqref{eq:rho} for $n \gg 1$), which is impossible. Consequently, $|\nabla \Psi_n(\cdot,t)|_2 < R_{0,n} \to R_0$ as $n \to +\infty$ from (\ref{g4}), so $\Psi_n$ is defined globally in time again from the well-posedness of \eqref{eq:time}. Therefore, $\bigl(\Psi_n(\cdot,t_n)\bigr)$ satisfies the assumptions of Lemma \ref{le:minim}, hence $\lim_n \inf_{u \in \cG} \|u - \Psi_n(\cdot,t_n)\| = 0$, a contradiction.
	\end{proof}
	
	From now on we assume (\ref{eq:os2}).

	Following \cite[Section 2.2]{Cazenave}, for $\phi \in \hrn$, we shall denote by $e^{\mathrm{i}t\Delta} \phi$ the unique strong solution\footnote{With a small abuse of notation, we shall write $e^{\mathrm{i}t\Delta} \phi$ for the map $t \mapsto e^{\mathrm{i}t\Delta} \phi$.} to the linear problem (which is global in time)
	\begin{equation*}
		\left\{ \begin{array}{l}
			\mathrm{i} \partial_t \Phi = -\Delta \Phi \\ 
			\Phi(\cdot,0) = \phi.
		\end{array} \right.
	\end{equation*}
	The solution $e^{\mathrm{i} t\Delta}\phi$ can be characterized using the Fourier variables, i.e.,
	$$
	e^{\mathrm{i} t\Delta}\phi (x) = \left( \frac{1}{4\pi \mathrm{i}t} \right)^{N/2} \int_{\R^N} e^{\frac{\mathrm{i}|x-y|^2}{4t}} \phi(y) \, \mathrm{d}y,
	$$
	and therefore
	$e^{\mathrm{i} t\Delta}\phi = K_t * \phi$, where 
	$$
	K_t (x) := \left( \frac{1}{4\pi \mathrm{i}t} \right)^{N/2} e^{\frac{\mathrm{i}|x|^2}{4t}}.
	$$
	
	Note that, under (\ref{F0}), the Nemytskii operator $N_f : H^1 (\R^N; \mathbb{C}) \rightarrow H^{-1} (\R^N; \mathbb{C})$
	$$
	u \mapsto N_f (u) := f \circ u
	$$
	is continuous and bounded over bounded subsets of $H^1 (\R^N; \mathbb{C})$, and a function $\Psi \in X(\infty,2,I)$ (respectively $\Psi \in \cC(I; H^1(\R^N; \mathbb{C}))$) is a weak (respectively strong) solution to \eqref{eq:time} if \textit{Duhamel's formula} hold. That is, from \cite[Proposition 3.1.3]{Cazenave}, $\Psi$ is a solution to \eqref{eq:time} if and only if
	\begin{equation*}
		\Psi(\cdot,t) = e^{\mathrm{i} t \Delta} \psi_0 + \mathrm{i}\int_0^t e^{\mathrm{i} (t-s) \Delta} f\bigl(\Psi(\cdot,s)\bigr) \, \ds
	\end{equation*}
	for (almost) every $t\ge0$.
	
	Now, inspired by \cite[Section 3]{JJLV}, we give the following technical results.

	Recalling that, if $p,r \in [2,+\infty]$, $(p,r)$ is an \emph{admissible pair} if and only if
	\[
	\frac{2}{p} + \frac{N}{r} = \frac{N}{2},
	\]
	with the very same proof we generalize \cite[Lemma 3.5]{JJLV} to a wider class of functions $f_\alpha$.
	\begin{Lem}\label{le:f_al}
		Let $\alpha \in(2,2^*]$ and define
		\[
		p_1 := \frac{4 \alpha}{(\alpha - 2) (N - 2)}, \quad r_1 := \frac{\alpha N}{\alpha + N - 2}.
		\]
		Then $(p_1,r_1)$ is an admissible pair and, moreover, for every admissible pair $(p_2,r_2)$ there exists $C>0$ such that for every $T>0$ and every $\Phi,\Psi \in X(p_1,r_1,T)$
		\[
		\left\| t \mapsto \int_0^t e^{\mathrm{i} (t - s) \Delta} \nabla [f_\alpha \circ \Phi(\cdot,s)] \, \ds \right\|_{Y(p_2,r_2,T)} \le C T^\mu \|\nabla \Phi\|_{Y(p_1,r_1,T)}^{\alpha-1}
		\]
		and
		\begin{multline*}
			\left\| t \mapsto \int_0^t e^{\mathrm{i} (t - s) \Delta} [f_\alpha \circ \Phi(\cdot,s) - f_\alpha \circ \Psi(\cdot,s)] \, \ds \right\|_{Y(p_2,r_2,T)}\\
			\le C T^\mu \left(\|\nabla \Phi\|_{Y(p_1,r_1,T)}^{\alpha-2} + \|\nabla \Psi\|_{Y(p_1,r_1,T)}^{\alpha-2}\right) \|\Phi(\cdot,s) - \Psi(\cdot,s)\|_{Y(p_1,r_1,T)},
		\end{multline*}
		where $\mu = \frac{(N-2)(2^*-\alpha)}{4} \ge 0$ and $f_\alpha \in W^{1,1}_\textup{loc}(\mathbb{R})$ is odd\footnote{This implies that $f_\alpha'$ is even, hence we define $f_\alpha'(z) := f_\alpha'(|z|)$ for $z \in \mathbb{C}$.} and satisfies $|f_\alpha'(t)| \lesssim |t|^{\alpha-2}$ for a.e. $t \in \mathbb{R}$.
	\end{Lem}
	\begin{Rem}
		We observe that, following the proof of \cite[Lemma 3.5]{JJLV}, the same result can be obtained for $\alpha = 2$ and $(p_1,r_1) = (+\infty,2)$.
	\end{Rem}
	
	From now on, we shall denote
	\begin{equation*}
		(p_q,r_q) := \left( \frac{4 q}{(q - 2) (N - 2)},\frac{q N}{\alpha + N - 2} \right) \quad \text{ and } \quad (\overline p,\overline r) := \left( 2^*,\frac{2 N^2}{N^2 - 2 N + 4} \right),
	\end{equation*}
	where $q$ is provided in Proposition \ref{pr:os} (observe that $(\overline p,\overline r)$ is given by $(p_q,r_q)$ replacing $q$ with $2^*$).\\
	For $T>0$, we define $Y(T) := Y(p_q,r_q,T) \cap Y(\overline p,\overline r,T)$ with norm $\|\cdot\|_{Y(T)} := \|\cdot\|_{Y(p_q,r_q,T)} + \|\cdot\|_{Y(\overline p,\overline r,T)}$ and $X(T) := X(p_q,r_q,T) \cap X(\overline p,\overline r,T)$ with norm $\|\cdot\|_{X(T)} := \|\cdot\|_{X(p_q,r_q,T)} + \|\cdot\|_{X(\overline p,\overline r,T)}$.
	
	The next lemma concerns, in particular, the local existence and uniqueness of a solution to \eqref{eq:time}. It is a preliminary step toward its global-in-time analogous. We point out that the presence of a generic $f$ makes the argument more involved.
	
	\begin{Lem}\label{le:33JJLV}
		There exists $\gamma>0$ such that if $\psi_0 \in H^1(\rn;\mathbb{C})$ and $T \in (0,1]$ satisfy $\|e^{\mathrm{i} t \Delta} \psi_0\|_{X(T)} \le \gamma$, then there exists a unique solution $\Psi \colon \rn \times [0,T] \to \mathbb{C}$ to \eqref{eq:time}. Moreover, $\Psi \in X(p,r,T)$ for every admissible pair $(p,r)$, and, for all $t \in [0,T]$, there holds $|\Psi(\cdot,t)|_2 = |\psi_0|_2$ and $J\bigl(\Psi(\cdot,t)\bigr) = J(\psi_0)$.
	\end{Lem}
	\begin{proof}
		Observe that (\ref{eq:os2}) implies that $|f(t) - f(s)| \lesssim |t - s| (|t|^{q-2} + |s|^{q-2} +|t|^{2^*-2}+|s|^{2^*-2})$ for all $t,s \in \R$.
		Let us write 
		$$
		f = f_1 + f_2
		\text{ with }
		f_j(t) := \int_0^t f_j'(s) \, \ds, \ (j = 1,2)
		\text{ and }
		f_1' = f' \chi_{\{|z| \le 1\}}.
		$$
		Since $f_j$ is odd and $f_j'$ is even, we extend them to the complex plane in the usual way. Then, for a.e. $z \in \mathbb{C}$, we have $|f_1'(z)| \lesssim |z|^{q-2}$ and $|f_2'(z)| \lesssim |z|^{2^*-2}$.\\
		For $\alpha,T>0$, let us denote by $B(\alpha,T)$ the closed ball of centre $0$ and radius $\alpha$ in $X(T)$. If we denote by $\mathfrak{d}$ the metric induced by $\|\cdot\|_{Y(T)}$, then $\bigl(B(\alpha,T),\mathfrak{d}\bigr)$ is a complete metric space in view of \cite[Lemma 3.7]{JJLV}.\\
		Let us also denote by $\mathcal{F}$ the operator that maps $\Phi \in X(T)$ to the function defined via
		\[
		\mathcal{F}(\Phi)(t) := e^{\mathrm{i} t \Delta} \phi + \mathrm{i} \int_0^t e^{\mathrm{i} (t-s) \Delta} [f \circ \Phi(\cdot,s)] \, \ds.
		\]
		From Lemma \ref{le:f_al} (used with $\Psi = 0$), for every $\Phi \in B(2\gamma,T)$ there holds
		\begin{multline*}
			\left\|t \mapsto \int_0^t e^{\mathrm{i} (t-s) \Delta} \nabla[f \circ \Phi(\cdot,s)] \, \ds\right\|_{Y(p,r,T)}\\
			\le \left\|t \mapsto \int_0^t e^{\mathrm{i} (t-s) \Delta} \nabla[f_1 \circ \Phi(\cdot,s)] \, \ds\right\|_{Y(p,r,T)} + \left\|t \mapsto \int_0^t e^{\mathrm{i} (t-s) \Delta} \nabla[f_2 \circ \Phi(\cdot,s)] \, \ds\right\|_{Y(p,r,T)}\\
			\lesssim \|\nabla \Phi\|_{Y(p_q,r_q,T)}^{q-1} + \|\nabla \Phi\|_{Y(\overline{p},\overline{r},T)}^{2^*-1} \lesssim \gamma^{q-1} + \gamma^{2^*-1} \lesssim \gamma^{q-1}
		\end{multline*}
		(we can assume $\gamma < 1$) and, analogously,
		\begin{align*}
			\left\|t \mapsto \int_0^t e^{\mathrm{i} (t-s) \Delta} [f \circ \Phi(\cdot,s)] \, \ds\right\|_{Y(p,r,T)} & \lesssim \|\nabla \Phi\|_{Y(p_q,r_q,T)}^{q-2} \|\Phi\|_{Y(p_q,r_q,T)} + \|\nabla \Phi\|_{Y(\overline{p},\overline{r},T)}^{2^*-2} \|\Phi\|_{Y(\overline{p},\overline{r},T)}\\
			& \lesssim \gamma^{q-1}.
		\end{align*}
		Taking $(p,r) = (p_q,r_q)$ and $(p,r) = (\overline{p},\overline{r})$ shows that, if $\gamma$ is sufficiently small, then
		\[
		\left\|\mathcal{F}(\Phi)\right\|_{X(T)} \le \gamma + C({N,q}) \gamma^{q-1} \le 2 \gamma,
		\]
		where $C({N,q}) > 0$ depends only on $N$ and $q$, i.e., $\mathcal{F}\bigl(B(2\gamma,T)\bigr) \subset B(2\gamma,T)$. Now, using again Lemma \ref{le:f_al} and the decomposition $f = f_1 + f_2$, for every $\Phi,\Psi \in B(2\gamma,T)$, there holds
		\begin{align*}
			\left\|\mathcal{F}(\Phi) - \mathcal{F}(\Psi)\right\|_{p,r} & \lesssim \left(\|\nabla \Phi\|_{p_q,r_q}^{q-2} + \|\nabla \Psi\|_{p_q,r_q}^{q-2}\right) \|\Phi - \Psi\|_{p_q,r_q} + \left(\|\nabla \Phi\|_{\overline{p},\overline{r}}^{2^*-2} + \|\nabla \Psi\|_{\overline{p},\overline{r}}^{2^*-2}\right) \|\Phi - \Psi\|_{\overline{p},\overline{r}}\\
			& \le 2 \left(\gamma^{q-2} \|\Phi - \Psi\|_{p_q,r_q} + \gamma^{2^*-2} \|\Phi - \Psi\|_{\overline{p},\overline{r}}\right) \lesssim \gamma^{q-2} \left(\|\Phi - \Psi\|_{p_q,r_q} + \|\Phi - \Psi\|_{\overline{p},\overline{r}}\right)
		\end{align*}
		As before, taking $(p,r) = (p_q,r_q)$ and $(p,r) = (\overline{p},\overline{r})$ shows that
		\[
		\left\|\mathcal{F}(\Phi) - \mathcal{F}(\Psi)\right\|_{Y(T)} \le C({N,q}) \gamma^{q-2} \|\Phi - \Psi\|_{Y(T)}
		\]
		and, if $\gamma$ is sufficiently small, $\mathcal{F}$ is a contraction in $\bigl(B(2\gamma,T),\mathfrak{d}\bigr)$. Then, from the Banach--Caccioppoli theorem, $\mathcal{F}$ has a unique fixed point in $B(2\gamma,T)$. The uniqueness in the whole space $X(T)$ is proved as in \cite[Proof of Proposition 3.3, Step 2]{JJLV}, while the property $\Phi \in \cC\bigl([0,T];H^1(\R^N; \mathbb{C})\bigr) \cap X(p,r,T)$ follows from the Strichartz estimates (cf. \cite[Theorem 2.3.3 and Remark 2.3.8]{Cazenave} and \cite[Corollary 1.4]{KeelTao}). Finally, the conservation laws follow from \cite[Propositions 1 and 2]{Ozawa}.
	\end{proof}
	
	Finally, arguing as in the proof of Proposition \ref{pr:os} under the assumption (\ref{eq:os1}), we have the following result, which, in fact, is even simpler because we do not need to justify the global existence.
	
	\begin{Lem}\label{le:osl}
		For every $\varepsilon > 0$ there exists $\delta > 0$ such that, if $\psi_0 \in \hrn$ satisfies $\inf_{u \in \cG} \|u - \psi_0\| \le \delta$, then $\sup_{0 \le t < T} \inf_{u \in \cG} \|u - \Psi(\cdot,t)\| \le \varepsilon$, where $\Psi$ is the unique solution to \eqref{eq:time} and $T>0$ is the maximal existence time of $\Psi$.
	\end{Lem}
	
	We can now conclude the proof of Proposition \ref{pr:os}.
	
	\begin{proof}[Proof of Proposition \ref{pr:os} if (\ref{eq:os2}) holds]
		Using Lemmas \ref{le:33JJLV} and \ref{le:osl}, we can argue as in \cite[Section 4]{JJLV}.
	\end{proof}
	
	\section{Positive-energy solution}\label{sec:pes}
	
	This section is devoted to the search for a positive-energy solution of \eqref{eq:main} and to the dynamics of the corresponding standing wave.
	
	As we did in Section \ref{sec:lm}, the assumptions in each of the following lemmas will be somewhat minimal. Thus, we start assuming from now on that
	\begin{equation}\label{eq:H-basic}
		H \text{ is of class $\cC^1$ and } h = H' \text{ satisfies } |h(t)| \lesssim |t| + |t|^{2^*-1} \text{ for every } t \in \R,
	\end{equation}
	and
	\[
	\text{there exists $\xi \ne 0$ such that $H(\xi) > 0$.}
	\]
	
	\subsection{The elliptic problem}
	
	First of all, let us observe that, if $u \in \hrn$ is such that $\displaystyle\int_{\R^N} H(u) \, \dx > 0$ (such functions exist, as shown in \cite[page 325]{BerLions}), then we can define
	\begin{equation}\label{eq:rDEF}
		r(u) := \left( \frac{N}{2  |\nabla u|_2^2} \int_{\rn} H(u) \, \dx\right)^{1/2},
	\end{equation}
	and so $u\bigl(r(u)\cdot\bigr) \in \cM$. This proves that $\cM \ne \emptyset$. In fact, arguing as in \cite[Lemma 4.1]{BMS} we easily obtain the following fact.
	
	\begin{Lem} 
		$\cM$ is a $\cC^1$-manifold of codimension 1 in $H^1 (\R^N)$.
	\end{Lem}
	
	In view of \eqref{phisec}, for every $u \in \cM$,
	\begin{equation}\label{eq:ineqInM-}
		\int_{\rn} 2_\# H(u) - h(u)u \, \dx < 0 \quad \mbox{if and only if} \quad u \in \cM_-.
	\end{equation}
	
	Now we give some preliminary technical results.
	
	\begin{Lem}\label{le:H1-H2}
		If $h_1(t)t \le a H_1(t)$ and $h_2(t)t \le 2^* H_2(t)$ for all $t \in \R$, then for every $u \in \cM_-$ there holds
		\begin{equation*}
			\frac{2_\# - a}{2^* - 2_\#} \int_{\rn} H_1 (u) \, \dx < \int_{\rn} H_2 (u) \, \dx.
		\end{equation*}
	\end{Lem}

	\begin{proof}
		Note that, from \eqref{eq:ineqInM-}, we have
		\begin{align*}
			(2_\# - a) \int_{\rn} H_1(u) \, \dx &\leq \int_{\rn} 2_\# H_1(u) - h_1(u)u \, \dx \\ &< \int_{\rn} h_2(u)u - 2_\# H_2 (u) \, \dx \leq (2^*-2_\#) \int_{\rn} H_2 (u) \, \dx. \qedhere
		\end{align*}
	\end{proof}

	\begin{Lem}\label{le:grad-positive}
		If $h_1(t)t \le a H_1(t)$, $h_2(t)t \le 2^* H_2(t)$, and $H_2(t) \lesssim |t|^b + |t|^{2^*}$ for every $t \in \R$ hold, then $\inf_{\cM_- \cap \cD_\rho} |\nabla u|_2 > 0$.
	\end{Lem}
	\begin{proof}
		Fix $u \in \cM_- \cap \cD_\rho$. Note that, from Lemma \ref{le:H1-H2}, we get
		$$
		|\nabla u|_2^2 = \frac{N}{2} \int_{\rn} H(u) \, \dx \lesssim \int_{\rn} H_2 (u) \, \dx \lesssim |u|_{b}^b + |u|_{2^*}^{2^*}.
		$$
		Therefore, from the Gagliardo--Nirenberg inequality,
		$$
		|\nabla u|_2^2
		\lesssim |\nabla u|_2^{\frac{N}{2}(b-2)} |u|_2^{b-\frac{N}{2}(b-2)} + |\nabla u|_2^{2^*}
		\lesssim |\nabla u|_2^{\frac{N}{2}(b-2)} + |\nabla u|_2^{2^*}.
		$$
		Note that $N(b-2)/2 > N(2_\# - 2)/2 = 2$ and therefore $|\nabla u|_2$ is bounded away from $0$.
	\end{proof}
	
	In some of the next lemmas, we will need the condition
	\begin{equation}\label{eq:condH}
		\lim_{|t| \to +\infty} \frac{H(t)}{|t|^{2^*}} = 0,
	\end{equation}
	which is a consequence, in particular, of (\ref{F3}), (\ref{F4}), and (\ref{H2n}).
	
	\begin{Lem}\label{le:Mempty}
		If (\ref{H0n}), (\ref{H2n}), and \eqref{eq:condH} hold and $\rho$ is sufficiently small, then $\cM_0 \cap \cD_\rho = \emptyset$.
	\end{Lem}
	\begin{proof}
		For $q \in (2,2^*)$, let us recall from \eqref{eq:GN} the definition of $C_{N,q}$.
		Suppose by contradiction that there is $u \in \cM_0 \cap \cD_\rho$, i.e., $u \in \cM \cap \cD_\rho$ and, by \eqref{phisec},
		\begin{equation}\label{eq:M0}
			\int_{\rn} 2_\# H(u) - h(u)u \, \dx = 0.
		\end{equation}
		Note that from \eqref{eq:M0} and (\ref{H2n}) we have
		\begin{align*}
			(2_\# - 2) \int_{\rn} H_1(u) \, \dx &\geq \int_{\rn} 2_\# H_1(u) - h_1(u)u \, \dx \\ &= \int_{\rn} h_2(u)u - 2_\# H_2 (u) \, \dx \geq (b-2_\#) \int_{\rn} H_2 (u) \, \dx.
		\end{align*}
		Hence, using that $u \in \cM$, we get
		\begin{align*}
			|\nabla u|_2^2
			&= \frac{N}{2} \int_{\rn} H(u) \, \dx
			\leq \frac{N}{2} \frac{b - 2}{b-2_\#} \int_{\rn} H_1(u) \, \dx
			\leq C \frac{N}{2} \frac{b - 2}{b-2_\#} \left( \rho^2 + |u|_a^a \right)
			\\
			&
			\leq C \frac{N}{2} \frac{b - 2}{b-2_\#} \left( \rho^2 + C_{N,a}^a \rho^{a - N(a-2)/2} |\nabla u|_2^{N(a-2)/2} \right) \\
			&= A_\rho + B_\rho |\nabla u|_2^{N(a-2)/2},
		\end{align*}
		where
		$$
		A_\rho := D \rho^2, \quad B_\rho := D C_{N,a}^a \rho^{a - N(a-2)/2}, \quad D := C \frac{N}{2} \frac{b - 2}{b-2_\#}.
		$$
		Observe that $A_\rho \to 0$ and $B_\rho \to 0$ as $\rho \to 0^+$. Therefore
		$$
		\frac{\sqrt{A_\rho} + \frac14}{\left( \sqrt{A_\rho} + \frac12 \right)^{N(a-2)/2}} \to \frac{1}{2^{2-N(a-2)/2}} > 0 \quad \mbox{as } \rho \to 0^+.
		$$
		Thus, for sufficiently small $\rho$ we have that
		$$B_\rho \leq \frac{\sqrt{A_\rho} + \frac14}{\left( \sqrt{A_\rho} + \frac12 \right)^{N(a-2)/2}}.$$
		Hence, Lemma \ref{lem:fromAbove} implies that
		\begin{align*}
			|\nabla u|_2 &\leq \sqrt{A_\rho} + \frac12 = \sqrt{D}\,\rho + \frac12.
		\end{align*}
		On the other hand,
		again from \eqref{eq:M0} and (\ref{H2n}), we have
		\begin{align*}
			(2^*-2_\#) \int_{\rn} H_2(u) \, \dx
			&\geq \int_{\rn} h_2(u)u - 2_\# H_2(u)  \, \dx \\
			&= \int_{\rn} 2_\# H_1 (u) - h_1(u)u  \, \dx
			\geq (2_\#-a) \int_{\rn} H_1 (u) \, \dx.
		\end{align*}
		Moreover, from (\ref{H0n}), (\ref{H2n}), and \eqref{eq:condH}, for every $\varepsilon > 0$ there exists $C_\varepsilon > 0$ such that
		\begin{equation*}
			H_2(t) \leq \varepsilon |t|^{2^*} + C_\varepsilon |t|^b,
		\end{equation*}
		and so,
		using (\ref{H2n}) again, we get
		\begin{align*}
			|\nabla u|_2^2 &= \frac{N}{2} \int_{\rn} H(u) \, \dx \leq \frac{N}{2} \frac{2^*-a}{2_\#-a} \int_{\rn} H_2 (u) \, \dx \\
			&\leq \frac{N}{2} \frac{2^*-a}{2_\#-a} \left( C_\varepsilon C_{N,b}^b |\nabla u|_2^{N(b-2)/2} \rho^{b-N(b-2)/2} + \varepsilon |\nabla u|_2^{2^*} S^{-2^*/2}\right).
		\end{align*}
		Using Lemma \ref{lem:fromBelow} we get that
		\begin{align*}
			|\nabla u|_2 \geq \min \left\{1, \left( \frac{N}{2} \frac{2^*-a}{2_\#-a} \left( C_\varepsilon C_{N,b}^b \rho^{b-N(b-2)/2} + \varepsilon S^{-2^*/2} \right) \right)^{\frac{2}{4-N(b-2)}} \right\}.
		\end{align*}
		Choose
		$$
		\varepsilon  \leq \frac{S^{2^*/2}}{N} \frac{2_\#-a}{2^*-a}
		$$
		and let $\rho$ be so small that
		$$
		\frac{N}{2} \frac{2^*-a}{2_\#-a} C_\varepsilon C_{N,b}^b \rho^{b-N(b-2)/2} \leq \frac12.
		$$
		Then 
		$$
		\left( \frac{N}{2} \frac{2^*-a}{2_\#-a} \left( C_\varepsilon C_{N,b}^b \rho^{b-N(b-2)/2} + \varepsilon S^{-2^*/2} \right) \right)^{\frac{2}{4-N(b-2)}} \geq 1
		$$
		because $4-N(b-2) < 0$. Thus,
		$$
		\sqrt{A_\rho} + \frac12 = \sqrt{D}\, \rho + \frac12 \geq |\nabla u|_2 \geq 1
		$$
		and we get a contradiction, by choosing sufficiently small $\rho$.
	\end{proof}
	
	In the next Lemmas, we show some properties of the functional $J$ restricted to $\cM_- \cap \cD_\rho$. We will need the condition
	\begin{equation}\label{eq:condF}
		\lim_{|t| \to +\infty} \frac{F(t)}{|t|^{2_\#}} = +\infty,
	\end{equation}
	which is a consequence, in particular, of (\ref{F4}) and (\ref{H2n}).
	
	\begin{Lem}\label{le:J-positive}
		If (\ref{F0}), (\ref{J1}), \eqref{eq:rho}, and \eqref{eq:condF} hold, then $\inf_{\cM_- \cap \cD_\rho} J > 0$.
	\end{Lem}
	\begin{proof}
		Let $u \in \cM_- \cap \cD_\rho$.
		From \eqref{eq:condF}, we have $\lim_{s \to +\infty} J(s\star u) = -\infty$, and so
		$1$, which is a local maximum point of $(0,+\infty) \ni s \mapsto J(s \star u) \in \R$ due to \eqref{equivM},
		is actually, from (\ref{J1}), the global maximum point.\\
		Now, let $s_{\max} > 0$ be the unique maximum point of $t \mapsto g(\rho,t) t^2$, guaranteed by \eqref{eq:rho}, where $g$ is defined by \eqref{eq:g}, see (\ref{max}) in Lemma \ref{le:g}, and choose $s_u > 0$ such that $|\nabla (s_u \star u)|_2 = s_{\max}$. Then, by \eqref{eq:Jg},
		$$
		J(u) = J(1 \star u) \geq J(s_u \star u) \geq g(\rho, |\nabla (s_u \star u)|_2) |\nabla (s_u \star u)|_2^2 = g(\rho, s_{\max}) s_{\max}^2 > 0,
		$$
		which completes the proof.
	\end{proof}
	
	\begin{Lem}\label{le:coercive}
		If (\ref{F0}), (\ref{F1}), (\ref{F3}), (\ref{J1}), \eqref{eq:rho}, and \eqref{eq:condF} hold, then $J|_{\cM_-^\textup{rad} \cap \cD_\rho}$ is coercive.
	\end{Lem}
	\begin{proof}
		Suppose that $(u_n) \subset \cM_-^\textup{rad} \cap \cD_\rho$ is a sequence such that $\|u_n\|\to +\infty$. Clearly, since $u_n \in \cD_\rho$, we get $|\nabla u_n|_2 \to +\infty$. We set $s_n := |\nabla u_n|_2^{-1}$ and note that $s_n \to 0^+$ as $n \to +\infty$. Define $v_n := s_n \star u_n$. Then $|v_n|_2=|u_n|_2$, so $v_n \in \cD_\rho$, and $|\nabla v_n|_2 = 1$. Hence, $(v_n)$ is bounded in $\hrnr$.
		Then, there exists $v \in H^1_\textup{rad}(\R^N)$ such that, up to a subsequence, $v_n \rightharpoonup v$ in $\hrnr$ and $v_n \to v$ in $L^{2_\#}(\R^N)$ and a.e. in $\R^N$. Let $F_\pm := \max\{\pm F,0\}$. From (\ref{F0}) and \eqref{eq:condF}, we obtain that $F_-(t) \lesssim t^2$. Suppose that $v\neq0$. Then, from Lemma \ref{le:J-positive}, \eqref{eq:condF}, and Fatou's lemma
		\begin{align*}
			0 & < \frac{J(u_n)}{|\nabla u_n|_2^2} = \frac12 - s_n^{N+2} \int_{\R^N} F(u_n(s_n x)) \, \dx = \frac12 + s_n^{N+2} \int_{\R^N} F_-(s_n^{-N/2} v_n) - F_+(s_n^{-N/2} v_n) \, \dx\\
			& \le \frac12 + C s_n^2 \rho^2 - \int_{\R^N} \frac{F_+(s_n^{-N/2}v_n)}{\bigl| s_n^{-N/2} v_n \bigr|^{2_\#}} |v_n|^{2_\#} \, \dx \to -\infty,
		\end{align*}
		which is a contradiction. Hence, $v=0$ and $v_n \to 0$ in $L^{2_\#}(\rn)$. Note that $u_n = s_n^{-1} \star v_n \in \cM_-^\textup{rad} \cap \cD_\rho$. Thus, from (\ref{J1}) and \eqref{eq:condF}, $s_n^{-1}$ is the unique maximizer of $s \mapsto J(s \star v_n)$. Consequently, we get
		$$
		J(u_n) = J(s_n^{-1} \star v_n)\geq J(s \star v_n) = \frac{s^2}{2} - s^{-N} \int_{\R^N} F (s^{N/2} v_n) \, \dx
		$$
		for any $s > 0$. Note that, from (\ref{F0}), (\ref{F1}), and (\ref{F3}), for every $\varepsilon > 0$ there is $C_\varepsilon > 0$ such that for every $n$
		$$
		\int_{\R^N} F (s^{N/2} v_n) \, \dx \leq \varepsilon (|s^{N/2} v_n|_2^2 + |s^{N/2} v_n|_{2^*}^{2^*}) + C_\varepsilon |s^{N/2} v_n|_{2_\#}^{2_\#}
		$$
		and therefore
		$$
		s^{-N} \int_{\R^N} F (s^{N/2} v_n) \, \dx \to 0 \quad \mbox{as } n\to +\infty,
		$$
		Hence
		$$
		\liminf_n J(u_n) \geq \liminf_n \left( \frac{s^2}{2} - s^{-N} \int_{\R^N} F (s^{N/2} v_n) \, \dx \right) = \frac{s^2}{2}
		$$
		for every $s > 0$. Thus $\lim_n J(u_n) = +\infty$.
	\end{proof}
	
	\begin{Rem}
		The statement of Lemma \ref{le:coercive} also holds if $\cM_-^\textup{rad}$ is replaced with $\cM_-$, i.e., $J|_{\cM_- \cap \cD_\rho}$ is coercive. The proof is similar, it is enough to replace the compact embedding of $H^1_{\mathrm{rad}} (\R^N)$ into $L^{2_\#} (\R^N)$ with Lions' concentration-compactness principle, see, e.g., \cite[Lemma 2.4]{BiegMed}.
	\end{Rem}
	
	In the next Lemma, we need the condition
	\begin{equation}\label{eq:supQ}
		\lim_{t \to 0} \frac{H(t)}{t^2} = 0,
	\end{equation}
	which is a consequence, in particular, of (\ref{F1}), (\ref{F4}), and (\ref{H2n}).
	
	\begin{Lem}\label{le:inf-achieved}
		If (\ref{F0}), (\ref{F1}), (\ref{F3}), (\ref{J1}), (\ref{H0n}), (\ref{H2n}), \eqref{eq:condH}, \eqref{eq:condF}, and \eqref{eq:supQ} hold and $\rho$ is sufficiently small, then $\inf_{\cM_-^\textup{rad} \cap \cD_\rho} J$ is attained; if $f$ is odd or $f|_{(-\infty,0)} \equiv 0$, then $\inf_{\cM_-^\textup{rad} \cap \cD_\rho} J = \inf_{\cM_- \cap \cD_\rho} J$ and it is attained by a non-negative and non-increasing (in the radial coordinate) function.
	\end{Lem}
	\begin{proof}
		Let $(u_n) \subset \cM_-^\textup{rad} \cap \cD_\rho$ be a minimizing sequence for $J$. From Lemma \ref{le:coercive}, $(u_n)$ is bounded in $H^1_{\mathrm{rad}}(\rn)$ and therefore, up to a subsequence, $u_n \rightharpoonup \tu$ in $H^1_{\mathrm{rad}}(\rn)$ and $u_n \to \tu$ in $L^p(\rn)$ for $2 < p < 2^*$ and a.e. in $\rn$ for some $\tu \in \cD_\rho$. This and $(u_n) \subset \cM$, together with (\ref{H0n}), \eqref{eq:condH}, and \eqref{eq:supQ}, imply that
		\begin{equation}\label{e:unot0}
			\frac{N}{2} \int_{\R^N} H(\tu) \, \dx = \lim_n \frac{N}{2} \int_{\R^N} H(u_n) \, \dx = \lim_n \int_{\R^N} |\nabla u_n|^2 \, \dx \ge \int_{\R^N} |\nabla \tu|^2 \, \dx.
		\end{equation}
		Additionally, $\tu \ne 0$ because, otherwise, \eqref{e:unot0} would yield $\lim_n |\nabla u_n|_2 = 0$, in contrast with Lemma \ref{le:grad-positive}.
		There follows
		\begin{equation*}
			0 < \int_{\rn} |\nabla \tu|^2 \, \dx \leq \frac{N}{2} \int_{\rn} H(\tu) \, \dx,
		\end{equation*}
		so we can define $\widetilde{r} := r(\tu) \ge 1$, cf. \eqref{eq:rDEF}, and $\tu \in \cM^\textup{rad}$. Note that
		$$
		\left| \tu (\widetilde{r} \cdot) \right|_2^2 = \widetilde{r}^{-N} | \tu|_2^2 \leq \widetilde{r}^{-N} \rho^2 \leq \rho^2,
		$$
		hence $\tu(\widetilde{r} \cdot) \in \cM^\textup{rad} \cap \cD_\rho$.\\
		Observe that, from (\ref{H2n}), \eqref{eq:condH} and \eqref{eq:supQ} still hold replacing $H(t)$ with $h(t)t$. Then, from \eqref{eq:ineqInM-},
		\begin{align*}
			\int_{\rn} 2_\# H\bigl(\tu (\widetilde{r} x)\bigr) - h\bigl(\tu (\widetilde{r} x)\bigr)\tu (\widetilde{r} x) \, \dx 
			& = \widetilde{r}^{-N} \int_{\rn} 2_\# H(\tu) - h(\tu)\tu \, \dx \\
			& = \widetilde{r}^{-N} \lim_n \int_{\R^N} 2_\# H(u_n) - h(u_n) u_n \, \dx \le 0,
		\end{align*}
		and so, from Lemma \ref{le:Mempty}, $\tu(\widetilde{r}\cdot) \in \cM_-^\textup{rad} \cap \cD_\rho$.\\ 
		In addition,
		\begin{align*}
			0
			< \inf_{\cM_-^\textup{rad} \cap \cD_\rho} J
			&\leq J\bigl(\tu (\widetilde{r} \cdot)\bigr)
			= \widetilde{r}^{-N} \int_{\rn}  \frac{N}{4} H(\tu)- F(\tu) \, \dx
			= \widetilde{r}^{-N} \lim_{n} J(u_n)
			\leq \lim_{n} J(u_n)
			= \inf_{\cM_-^\textup{rad} \cap \cD_\rho} J,
		\end{align*}
		hence $\widetilde r = 1$, and $\tu \in \cM_-^\textup{rad} \cap \cD_\rho$, and $J(\tu) = \inf_{\cM_-^\textup{rad} \cap \cD_\rho} J$.\\
		Now, assume that $f$ is odd\footnote{Let us recall that, in this case, we are considering complex-valued functions.}. If we denote by $u^*$ the Schwarz rearrangement of $|\widetilde u|$ (cf. \cite[Chapter 3]{LiebLoss}), we have that
		\[
		|u^*|_2 = |\widetilde u|_2,
		\quad
		|\nabla u^*|_2 \le |\nabla \widetilde u|_2,
		\]
		\[
		\int_{\rn} F(u^*) \, \dx = \int_{\rn} F(\widetilde u) \, \dx,
		\quad
		\int_{\rn} H(u^*) \, \dx = \int_{\rn} H(\widetilde u) \, \dx,
		\quad
		\int_{\rn} h(u^*)u^* \, \dx = \int_{\rn} h(\widetilde u) \widetilde u \, \dx,
		\]
		which implies that $r(u^*) \ge 1$, see \eqref{eq:rDEF}. Then, arguing as above, we obtain that $r(u^*) = 1$. Therefore, using also \eqref{eq:ineqInM-}, $u^* \in \cM_-^\textup{rad}$, $|\nabla u^*|_2 = |\nabla \widetilde u|_2$, and $J(u^*) = J(\widetilde u)$. If $f|_{(-\infty,0)} \equiv 0$, we consider the Schwarz rearrangement of $\max\{\tu,0\}$ and then an almost identical argument applies.
	\end{proof}
	
	Now we can prove the existence of a second solution to \eqref{eq:main}.
	\begin{proof}[Proof of Theorem \ref{th:2sol}]
		Let $\tu \in \cM_-^\textup{rad} \cap \cD_\rho$ be the minimizer found in Lemma \ref{le:inf-achieved}.
		Arguing as in \cite[Proof of Lemma 5.2]{Soave} and using Lemma \ref{le:Mempty}, it can be proved that for every $v \in \cS \cap \cM^\textup{rad}$ the functional $\bigl(\Phi'(v),M'(v)\bigr) \colon \hrnr \to \R^2$ is surjective, where $\Phi(v) := |v|_2^2 - \rho^2$ and $M$ is defined in \eqref{defM}.
		Then, from \cite[Proposition A.1]{MedSch}, there exist $\lambda_{\tu} \ge 0$ and $\mu_{\tu} \in \R$
		such that $\tu$ is a weak solution to
		\begin{equation*}
			-(1+\mu_{\tu}) \Delta \tu + \lambda_{\tu} \tu = f(\tu) + \mu_{\tu} \frac{N}{4} h(\tu).
		\end{equation*}
		Hence,
		\begin{equation}\label{eq:nehari}
			(1+\mu_{\tu}) |\nabla \tu|_2^2 + \lambda_{\tu} |\tu|_2^2 = \int_{\rn} f(\tu)\tu + \mu \frac{N}{4} h(\tu)\tu \, \dx.
		\end{equation}
		If $\mu_{\tu} = -1$, then from \eqref{eq:nehari} and $\tu \in \cM_-$ we get
		\begin{align*}
			0
			&\le \lambda_{\tu} |\tu|_2^2
			= \int_{\rn} f(\tu)\tu - \frac{N}{4} h(\tu)\tu \, \dx
			< \int_{\rn} f(\tu)\tu - \frac{N}{4} 2_\# H(\tu) \, \dx \\
			&= -\frac{N}{2} \int_{\rn} H(\tu) \, \dx
			+ \int_{\rn} f(\tu)\tu - H(\tu) \, \dx
			= -2J(\tu) < 0,
		\end{align*}
		which is a contradiction. Hence, $\mu_{\tu} \neq -1$ and $\tu$ satisfies the Poho\v{z}aev identity
		\begin{equation}\label{poho}
			(1+\mu_{\tu}) |\nabla \tu|_2^2 = 2^* \int_{\R^N} F(\tu) + \mu_{\tu} \frac{N}{4} H(\tu) - \frac{\lambda_{\tu}}{2} |\tu|^2 \, \dx.
		\end{equation}
		Combining $\tu \in \cM_-$, \eqref{eq:nehari}, and \eqref{poho}, we have
		\begin{equation*}
			(1+\mu_{\tu}) \frac{N}{2} \int_{\rn} H(\tu) \, \dx
			= (1+\mu_{\tu}) |\nabla \tu|_2^2
			= \frac{N}{2} \int_{\rn} H(\tu) + \frac{N}{4} \mu_{\tu} \bigl(h(\tu)\tu - 2 H(\tu)\bigr) \, \dx
		\end{equation*}
		or, equivalently,
		$$
		\mu_{\tu} \int_{\rn} 2_\# H(\tu) - h(\tu)\tu \, \dx = 0.
		$$
		Hence, \eqref{eq:ineqInM-} implies that $\mu_{\tu} = 0$.\\
		Therefore, $\tu$ solves
		$$
		-\Delta \tu + \lambda_{\tu} \tu = f(\tu).
		$$
		If $|\tu|_2 < \rho$, then $\lambda_{\tu} = 0$ and $\tu$ is a weak solution to 
		$$
		-\Delta \tu = f(\tu).
		$$
		Thus,
		\[
		\int_{\rn} f(\tu)\tu \, \dx
		= |\nabla \tu|_2^2
		= \frac{N}{2} \int_{\rn} H(\tu) \, \dx
		\]
		or, equivalently,
		$$
		\int_{\rn} 2^* F(\tu) - f(\tu)\tu \, \dx = 0.
		$$
		From (\ref{F4}),
		$$
		2^* F(\tu(x)) - f(\tu(x))\tu(x) = 0 \quad \mbox{for a.e. } x \in \R^N,
		$$
		while, from Remark \ref{rem:12} (\ref{Linfty}), $\tu$ is continuous and $u(x) \to 0$ as $|x|\to+\infty$.
		Hence, there exists an open interval $I \subset \R$ such that $0 \in \overline{I}$ and for all $t \in I$
		$$
		2^* F(t) - f(t)t = 0.
		$$
		Therefore, $F(t) = C |t|^{2^*}$ for $t \in I$, which is a contradiction with (\ref{F2}). Hence, $\lambda_{\tu} > 0$ and $\tu \in \cS$ is a solution to
		\[
		-\Delta \tu + \lambda_{\tu} \tu = f(\tu).
		\]
		If $f$ is odd (respectively, $f|_{(-\infty,0)} \equiv 0$), then from Lemma \ref{le:inf-achieved} we can replace $\widetilde u$ with the Schwarz rearrangement of its modulus (respectively, its positive part), which is non-negative and non-increasing in the radial coordinate. The strong minimum principle and (\ref{F4}) yield that $\widetilde u$ is positive. As for the final part, assume $\widetilde u$ is constant in the annulus $\mathcal{A} := \{r_1 \le |x| \le r_2\}$ for some $r_2 > r_1 > 0$. Then, $f(\widetilde u) - \lambda_{\tu} \widetilde u = -\Delta \widetilde u = 0$ in $\mathcal{A}$ and, since $\widetilde u$ is non-increasing, $-\Delta \widetilde u = f(\widetilde u) - \lambda_{\tu} \widetilde u \le 0$ in $\{|x| \ge r_1\}$. At the same time, $\widetilde u$ attains its maximum at every interior point of $\mathcal{A}$, which contradicts the strong maximum principle.
	\end{proof}
	
	\subsection{Dynamics}
	
	Let us give sufficient conditions for finite-time blow-up to occur.
	
	\begin{Lem}\label{le:inst}
		Under the assumptions of Proposition \ref{pr:si}, if $u \in \cD_\rho \setminus \{0\}$ satisfies $J(u) < \inf_{\cM_- \cap \cD_\rho} J$, $t_u < 1$, and $|\cdot|u \in L^2(\rn)$, then $u$ is strongly unstable.
	\end{Lem}
	\begin{proof}
		Let $u \in \cD_\rho \setminus \{0\}$ be as in the assumptions.
		Since for $s>0$
		\[
		\frac{\mathrm{d}}{\mathrm{d}s} J(s\star u) = \frac{M(s\star u)}{s},
		\]
		from (\ref{J2}) we have
		\begin{equation*}
			M(u) = M(1\star u) < \frac{M(t_u \star u)}{t_u} = 0
		\end{equation*}
		and so, since
		\[
		J(1\star u) - J(t_u \star u)
		=
		\frac{M(\xi\star u)}{\xi} (1-t_u)
		\geq
		M(u) (1-t_u)
		\]
		with $\xi\in[t_u,1]$,
		\begin{equation}\label{eq:JinfJM}
			J(u) = J(1\star u)
			\ge \inf_{\cM_- \cap \cD_\rho} J + M(u) (1 - t_u)
			\ge \inf_{\cM_- \cap \cD_\rho} J + M(u).
		\end{equation}
		Now let us consider \eqref{eq:time} with $\psi_0 = u$ and its solution $\Psi$ (to lighten the notations, we will write $\Psi(s)$ instead of $\Psi(\cdot,s)$).\\
		The continuity of $\cD_\rho\setminus\{0\} \ni v \mapsto t_v \in (0,+\infty)$ (applying the implicit function theorem and using that $t_v\star v\in \cM_-\cap \cD_\rho$) implies that there exists $\delta>0$ such that $M(\Psi(s)) \le -\delta$ for every $s \in [0,T)$, where $T \in (0,+\infty]$ is the maximal existence time. Indeed, since, by assumptions, $t_{\Psi(0)} = t_u< 1$, assume by contradiction that there exists $\bar s>0$ such that $t_{\Psi(s)} < 1$ for all $s\in[0,\bar s)$ and $t_{\Psi(\bar s)} = 1$. Then $M(\Psi(\bar s)) = 0$. On the other hand, from \eqref{eq:JinfJM} and the conservation of mass and energy, there holds
		\begin{equation}\label{Mpsis}
			M\bigl(\Psi(s)\bigr) \le J\bigl(\Psi(s)\bigr) - \inf_{\cM_- \cap \cD_\rho} J = J(u) - \inf_{\cM_- \cap \cD_\rho} J =: -\delta < 0
		\end{equation}
		for every $s \in [0,\bar{s})$, and by continuity, $M(\Psi(\bar s)) \le -\delta$ as well, reaching a contradiction.\\
		Thus, \eqref{Mpsis} applies to every $s \in [0,T)$.\\
		From the virial identity (cf. \cite[Proposition 6.5.1]{Cazenave}), the function
		\[
		s \mapsto V(s) := \int_{\R^N} |x|^2 |\Psi(x,s)|^2 \, \dx
		\]
		is of class $\cC^2$ and $V''(s) = 8M\bigl(\Psi(s)\bigl) \le -8 \delta$ for every $s \in [0,T)$. In particular
		$$0 \le V(s) \le V(0) + V'(0) s - 4 \delta s^2 \text{ for every } s \in [0,T).$$
		Since the right-hand side diverges negatively as $s \to +\infty$, we have $T < +\infty$, which implies that $\lim_{s \to T^-} |\nabla \Psi(s)|_2 = +\infty$.
	\end{proof}
	
	We can now prove the last result of this work.
	
	\begin{proof}[Proof of Proposition \ref{pr:si}]
		For $s>1$, let $\Psi_s$ be the solution to \eqref{eq:time} with $\psi_0 = s\star \widetilde{u}$. Since $s\star \widetilde{u} \to \widetilde u$ in $\hrn$ as $s \to 1^+$, it suffices to show that $|\nabla\Psi_s|_2$ blows up in finite time for every $s>1$.
		From (\ref{J1}), the unique maximizer $t_{s\star \widetilde{u}}$ of $t\mapsto J(t\star(s\star\widetilde{u}))$
		satisfies $st_{s\star \widetilde{u}}=1$ since $J(t\star(s\star\widetilde{u}))=J((ts)\star \widetilde{u})$ and $\widetilde{u}\in \mathcal{M}_-$. Thus,
		$t_{s\star \widetilde{u}}=1/s<1$ and
		\[
		J(s\star \widetilde u) < J(\widetilde u) = \inf_{\cM_- \cap \cD_\rho} J.
		\]
		Moreover, since $\lambda_{\widetilde{u}} > 0$ and $\widetilde u$ is radial, from \cite[Lemma 2]{BerLions} we have that $\widetilde u$ decays exponentially at infinity, hence $|\cdot| \widetilde u \in L^2(\rn)$. The assertion then follows from Lemma \ref{le:inst}.
	\end{proof}
	
	\appendix
	
	\section{Useful estimates}
	In this Appendix, we prove the estimates that played an important role in the crucial Lemma \ref{le:Mempty}.
	\begin{Lem}\label{lem:fromAbove}
		Suppose that $x > 0$ satisfies $x^2 \leq A + B x^{p}$, where $p \in (0,2)$ and $A, B > 0$ satisfy $B( \sqrt{A} + 1/2 )^{p} \leq \sqrt{A} + 1/4 $. Then
		$$
		x \leq \sqrt{A} + \frac12.
		$$
	\end{Lem}
	
	\begin{proof}
		We introduce the function $f \colon [0, +\infty) \rightarrow \R$ by
		$$
		f(t) := t^2 - A - B t^p.
		$$
		Note that $f$ is of class $\cC^1$ on $(0,+\infty)$ and
		$$
		f'(t) = 2t - Bp t^{p-1} = t \Big( 2-\frac{pB}{t^{2-p}} \Big).
		$$
		Hence, $f$ is decreasing on $[0, t_0]$ and increasing to $+\infty$ on $[t_0, +\infty)$, with $t_0:=(pB/2)^{1/(2-p)}$. Consequently, there exists a unique $t_1 > t_0$ such that $f$ is negative on $(0,t_1)$ and positive on $(t_1, +\infty)$. 
		Moreover, observe that
		$$
		f \left(\sqrt{A} + \frac12 \right) \geq \sqrt{A} + \frac14 - \sqrt{A} - \frac14 = 0.
		$$
		Since $f(x) \leq 0$, we conclude.
	\end{proof}
	
	\begin{Lem}\label{lem:fromBelow}
		Suppose that $x > 0$ satisfies $x^2 \leq A x^{p} + B x^q$, where $2 < p < q$ and $A, B > 0$. Then
		$$
		x \geq \min \left\{ 1, (A+B)^{1/(2-p)} \right\}.
		$$
	\end{Lem}
	
	\begin{proof}
		Let us rewrite the given inequality as
		$$
		1 \leq A x^{p-2} + B x^{q-2}
		$$
		and consider the function $f \colon [0, +\infty) \rightarrow \R$ given by
		$$
		f(t) = 1 - A t^{p-2} - B t^{q-2}.
		$$
		If we compute 
		$$
		f'(t) = -A (p-2) t^{p-3} - B (q-2) t^{q-3} = - t^{p-3} ( A(p-2) + B(q-2) t^{q-p}),
		$$
		it is clear that $f$ is decreasing. Moreover, $f(x) \leq 0$. Denote 
		$$
		\xi := \min \left\{ 1, (A+B)^{1/(2-p)} \right\}.
		$$
		If $A+B < 1$, then $\xi = 1$ and
		$$
		f(\xi) = f(1) = 1-A-B > 0,
		$$
		so $x \geq \xi$.\\
		On the other hand, if $A+B \geq 1$, then $\xi = (A+B)^{1/(2-p)} \leq 1$ and
		$$
		f(\xi) = 1 - \frac{A}{A+B} - B \left( \frac{1}{A+B} \right)^{(q-2)/(p-2)} \geq 1 - \frac{A}{A+B} - \frac{B}{A+B} = 0,
		$$
		so $x \geq \xi$.
	\end{proof}

	\section{Examples for (\ref{J1}) and (\ref{J2}).}\label{Example}
	In the first part of this appendix, we provide an example of a nonlinear term $f$ that satisfies (\ref{J1}) and, under additional assumptions,  (\ref{J2}). It generalizes the case given by two different powers.\\
	In the second part, we present some sufficient conditions on the nonlinear term for (\ref{J2}) to hold, paired with an example that does not consist merely of powers.
	\subsection{Multiple powers}\label{appB1}
	
	Consider the nonlinearity $f \colon \R \rightarrow \R$ given by
	$$
	f(t) := \sum_{k=0}^K \alpha_k |t|^{q_k - 2} t
	+ \sum_{\ell=0}^L \beta_\ell |t|^{p_\ell - 2} t,
	$$
	where $2 < q_0 < \dots < q_K < 2_\# < p_0 < \dots < p_L \le 2^*$, $q_k, p_\ell \in \mathbb{Q}$, $\alpha_k, \beta_\ell > 0$, and $K, L \geq 0$. Fix $u \in \cD_\rho \setminus \{0\}$, with $\rho$ satisfying \eqref{eq:rho}. Note that
	$$
	\varphi(s) := J(s \star u) = \frac{|\nabla u|_2^2}{2} s^2 - \sum_{k=0}^K \frac{\alpha_k}{q_k}  s^{N (q_k-2) / 2}|u|_{q_k}^{q_k} - \sum_{\ell=0}^L \frac{\beta_\ell}{p_\ell} s^{N (p_\ell-2) / 2}  |u|_{p_\ell}^{p_\ell}
	$$
	and
	\begin{align*}
		\varphi'(s)
		&=
		|\nabla u|_2^2 s
		- \sum_{k=0}^K \frac{N}{2} (q_k-2) \frac{\alpha_k}{q_k} s^{N (q_k-2) / 2 - 1}|u|_{q_k}^{q_k}
		- \sum_{\ell=0}^L\frac{N}{2} (p_\ell-2) \frac{\beta_\ell}{p_\ell} s^{N (p_\ell-2) / 2 - 1} |u|_{p_\ell}^{p_\ell} \\
		&= s^{-1}
		\left(
		- \sum_{k=0}^K \frac{N}{2} (q_k-2) \frac{\alpha_k}{q_k} s^{N (q_k-2) / 2} |u|_{q_k}^{q_k}
		+ |\nabla u|_2^2 s^2
		- \sum_{\ell=0}^L \frac{N}{2} (p_\ell-2) \frac{\beta_\ell}{p_\ell} s^{N (p_\ell-2) / 2} |u|_{p_\ell}^{p_\ell}
		\right).
	\end{align*}
	Note also that
	$$
	0 < \frac{N}{2} (q_k-2) < 2 < \frac{N}{2} (p_\ell-2),
	$$
	thus, from Lemma \ref{le:g}, $\varphi$ has a global maximum point at a positive level; hence, it also has a local minimum point at a negative level.\\
	Let us write
	$$\frac{N}{2} (q_k-2) = \frac{a_k}{b_k},
	\quad
	\frac{N}{2} (p_\ell-2) = \frac{c_\ell}{d_\ell}
	$$
	for some $a_k, b_k, c_\ell, d_\ell \in \mathbb{N}$ such that $\operatorname{gcd} \left( a_k, b_k \right) = \operatorname{gcd} \left( c_\ell, d_\ell \right) = 1$, and $m := \operatorname{lcm} \left( b_0, \ldots, b_K, d_0, \ldots, d_L \right)$. Then,
	$$
	\varphi'(s) = s^{-1} P(s^{1/m}),
	$$
	where $P$ is the polynomial given by
	$$
	P(t)
	:=
	- \sum_{k=0}^K \frac{\alpha_k}{q_k} \frac{a_k}{b_k} t^{ m a_k / b_k}|u|_{q_k}^{q_k}
	+ |\nabla u|_2^2 t^{2m}
	- \sum_{\ell=0}^L \frac{\beta_\ell}{p_\ell} \frac{c_\ell}{d_\ell} t^{m c_\ell / d_\ell} |u|_{p_\ell}^{p_\ell}.
	$$
	From Descartes' rule of signs, $P$ has at most two positive roots. Thus, from the argument above, $\varphi'$ has exactly two roots, one of them is a local minimum point of $\varphi$ and one of them, $t_u$, is a local (hence global) maximum point of $\varphi$. Therefore, (\ref{J1}) is satisfied.\\
	Now, to see whether this nonlinearity satisfies (\ref{J2}), let us compute
	\[
	\varphi''(s) = s^{-2} \Biggl(-\sum_{k=0}^K \frac{\alpha_k}{q_k} \frac{a_k}{b_k} \biggl(\frac{a_k}{b_k} -1\biggr) |u|_{q_k}^{q_k} s^{a_k/b_k}
	+ |\nabla u|_2^2 s^2 - \sum_{\ell=0}^L \frac{\beta_\ell}{p_\ell} \frac{c_\ell}{d_\ell} \biggl(\frac{c_\ell}{d_\ell} - 1\biggr) |u|_{p_\ell}^{p_\ell} s^{c_\ell/d_\ell}\Biggr)
	\]
	and observe that, for a real number $r$, $N(r-2)/2-1 \lesseqgtr 0$ if and only if $r \lesseqgtr 2 + 2/N$.
	Then $\varphi$ is convex (respectively, concave) in a right-hand neighbourhood of the origin if $q_0 < 2+2/N$ (respectively, $q_0 > 2+2/N$), and, if $q_0 = 2 + 2/N$, then $\varphi$ is convex (respectively, concave) in a right-hand neighbourhood of the origin if $K=0$ (respectively, $K\ge1$).\\
	Let us also recall that $\varphi$ is convex in a neighbourhood of its local minimizer, concave in a neighbourhood of its global maximizer, and that, therefore, its second derivative changes sign between such two points.\\
	Now, let us consider the following cases for the sign pattern of the non-zero coefficients of $\varphi''$, i.e. the sequence of the signs of the coefficients ordered by ascending variable exponent. \\
	If $q_K \le 2+2/N$, then the sign pattern is $+,\dots,+,-,\dots,-$ and $\varphi''$ changes sign once on $(0,t_u)$, hence $\varphi$ is concave on $(t_u,+\infty)$. In fact, in this case, as we will see in Appendix \ref{appB2} below, the assumption that $q_k$'s and $p_\ell$'s are rational can be dropped.\\
	If $q_0 > 2+2/N$ or $q_0 = 2+2/N$ and $K\geq 1$, then the sign pattern is $-,\dots,-,+,-,\dots,-$ and $\varphi''$ changes sign twice on $(0,t_u)$, hence $\varphi$ is concave on $(t_u,+\infty)$.\\
	Finally, a similar argument applies to the case when the powers $q_k$ ($0 \le k \le K$) and $p_\ell$ ($0 \le \ell \le L$) are (positive and) rational multiples of a given real number.
	
	\subsection{A logarithmic term and concrete assumptions}\label{appB2}
	We introduce 
	$$
	G(t) :=  h(t)t - \left(2+\frac{2}{N} \right) H(t), \quad t \in \R,
	$$
	where $H$ is given by \eqref{eq:H} and satisfies \eqref{eq:H-basic}.
	If $G$ satisfies
	\begin{enumerate}[label=(G\arabic{*}),ref=G\arabic{*}]\setcounter{enumi}{-1}
		\item \label{G0} $G$ is even,
		\item \label{G1} $\limsup_{t \to 0} \frac{G(t)}{|t|^{2_\#}} \leq 0$,
		\item \label{G2} $\lim_{|t| \to +\infty} \frac{G(t)}{|t|^{2_\#}} = +\infty$,
		\item \label{G3} $t \mapsto \frac{G(t)}{t^{2_\#}}$ is increasing on $(0,+\infty)$,
	\end{enumerate}
	then, for fixed $u \in H^1(\R^N) \setminus\{0\}$,
	$$
	\frac{\mathrm{d}^2}{\mathrm{d} s^2} J(s \star u) = |\nabla u|_2^2 - \frac{N^2}{4} s^{-N-2} \int_{\R^N} G(s^{N/2} u) \, \dx = |\nabla u|_2^2 - \frac{N^2}{4} \int_{\R^N} \frac{G(s^{N/2} u)}{(s^{N/2})^{2_\#}} \, \dx
	$$
	has exactly one positive zero, and (\ref{J2}) holds provided $\rho$ satisfies \eqref{eq:rho}.
	
	In this case, we propose the following example
	$$
	F(t) = \frac{3}{7} |t|^{7/3} \ln (e+|t|)  + \frac{3}{13} |t|^{13/3}, \quad N=3, \ 2 < \frac73 < 2_\# = \frac{10}{3} < \frac{13}{3} <  2^*= 6,
	$$
	which satisfies (\ref{F0})--(\ref{F4}). Then
	$$
	G(t)
	=
	- \frac{1}{21} |t|^{7/3} \ln (e+|t|)
	+ \frac{3}{7} \frac{|t|^{10/3}}{e+|t|}
	- \frac{3}{7} \frac{|t|^{13/3}}{(e+|t|)^2}
	+ \frac{35}{39} |t|^{13/3}
	$$
	satisfies (\ref{G0})--(\ref{G3}).
	
	Finally, we observe that, if we consider $f$ as in Appendix \ref{appB1} with $q_K \le 2+2/N$, then $f$ satisfies (\ref{G0})--(\ref{G3}) even for real exponents $q_k$, $p_\ell$.

	\addtocontents{toc}{\protect\setcounter{tocdepth}{1}}
	
	\subsection*{Acknowledgements}
	P. d'Avenia and J. Schino are members of GNAMPA (INdAM) and are supported by the GNAMPA project {\em Metodi variazionali e topologici per alcune equazioni di Schr\"odinger nonlineari}. 
	They were also supported by the National Science Centre, Poland (grant no. 2020/37/N/ST1/00795).
	P. d'Avenia was also supported by PRIN 2017JPCAPN {\em Qualitative and quantitative aspects of nonlinear PDEs} and by European Union - Next Generation EU - PRIN 2022 PNRR P2022YFAJH {\em Linear and Nonlinear PDE’s: New directions and Applications}.
	B. Bieganowski was supported by the programme \textit{Excellence Initiative - Research University} at the University of Warsaw, New Ideas 3A (grant no. 01/IDUB/2019/94). Moreover, part of this work was completed during short research stays at Dipartimento di Meccanica, Matematica e Management, Politecnico di Bari and Institute of Mathematics, Polish Academy of Sciences; the authors wish to thank both institutions for the warm hospitality.

\end{document}